\documentclass[11pt,a4paper]{article}
\usepackage[T1]{fontenc}
\usepackage[latin1]{inputenc}
\usepackage[french]{babel}
\DeclareMathSizes{11}{11}{7}{5}
\usepackage{amsfonts}
\usepackage{epsfig}
\usepackage{bm}
\usepackage{bbm}
\usepackage{amsmath,amsthm,amssymb}
\usepackage{tabularx}
\usepackage{multirow}
\include{graphics}
\usepackage{mathrsfs}
\usepackage{natbib} 
 \usepackage{fancyhdr}

\def\be{\begin{equation} \displaystyle}
\def\ee{\end{equation} }

\def\bea{\begin{eqnarray}}
\def\eea{\end{eqnarray} }
\def\bean{\begin{eqnarray*}}
\def\eean{\end{eqnarray*} }

\def\CC{\rm \hbox{C\kern-.57em\raise.47ex  
         \hbox{$\scriptscriptstyle |$}\kern+0.5 em }}
\def\s{\scriptscriptstyle}

\def\bl{\bigl}
\def\br{\bigr}
\def\Bl{\Bigl}
\def\Br{\Bigr}
\DeclareMathOperator{\tr}{tr}

\theoremstyle{plain}
\newtheorem{theoreme}{Théorème}[section]
\newtheorem{corollaire}{Corollaire}[section]
\newtheorem{lemme}{Lemme}[section]
\newtheorem{definition}{Définition}[section]
\newtheorem{proposition}{Proposition}[section]
\newtheorem{remarque}{Remarques}
\newtheorem{remark}{Remarque}
\title{{\textbf {Théorèmes Limites Avec Poids Pour Les Martingales Vectorielles à Temps Continu  } }}
\author{\textsc{\textbf{  Faouzi Chaabane\footnote{\'Equipe d'Analyse Stochastique et Mod\'elisation Statistique (DGRST, E07/C15) Faculté Des Sciences de Bizerte. 7021 Jarzouna, Tunisie.}$\;$ 
 \& Ahmed Kebaier\footnote{Laboratoire d'Analyse et de Mathématiques Appliquées,  UMR 8050, Université  de Marne-La-Vallée, 5 boulevard Descartes,  Champs-Sur-Marne, F-77454 Marne-La-Vallée Cedex 2, France.}}}}

\begin{document}

%
%
%
%
%
%
%

\maketitle
%



\noindent \textbf{R\'{e}sum\'{e}}: On développe une approche générale du Théorème limite centrale presque sure 
 pour les martingales vectorielles quasi-continues à gauches convenablement normalisées et on dégage une extension 
quadratique de ce théorème tout en précisant les vitesses de convergence qui lui sont associés.L'application de 
ce résultat à un P.A.I.S. illustre l'usage qu'on peut en faire en statistique.  
\section{Introduction.}

\subsection{Motivation.}

\'Etabli suite aux  travaux pionniers de \cite{bros} et de \cite{schatte}, le théorème de la limite centrale presque-sûre (TLCPS)
a révélé
un nouveau phénomène dans la théorie classique des théorèmes limites. En effet, pour une  marche aléatoire $(S_n)_{n\geq 1}$ 
à valeurs dans $\mathbb R^d$ et dont les accroissements sont des v.a. i.i.d., centrés
de variance $C$, le (TLCPS) assure   que  les mesures empiriques logarithmiques associées aux v.a. $(n^{-1/2}S_n)$ c'est à dire:
$$\mu_N=(\log N)^{-1}\sum_{n=1}^N n^{-1}\delta_{n^{-1/2}S_n}$$
vérifient
$$\mu_N\Rightarrow \mu_{\infty}\;\;p.s.,$$
o\`u $\mu_{\infty}$ est la loi Gaussienne de moyenne $0$, de variance $C$ et $\delta_x$ 
la mesure de Dirac en $x$. Dans ce cadre, et sous des conditions d'uniformes intégrabilité par 
exemple,
on dispose de la propriété suivante appelée loi forte quadratique (LFQ):
$$\lim_{N\rightarrow \infty} (\log N)^{-1} \sum_{n=1}^N  n^{-2}S_nS_n^*=C\;\;p.s.,$$
o\`u $S_n^*$ désigne le transposé du vecteur $S_n$.

Le  théorème de la limite centrale presque-sûre ainsi que les divers théorèmes ``logarithmiques'' qui lui sont associés ont  mené à 
une littérature étendue durant  
la décennie passée. En effet, ils  ont été généralisés aux martingales 
discrètes unidimensionnelles par \cite{faouzi1} et \cite{lifshits} puis 
aux martingales discrètes $d$-dimensionnelles par \cite{f&f&a} et  ensuite aux martingales continues par \cite{faouzic}.

Les résultats de \cite{faouzi1} et de \cite{faouzic} ont été obtenus gr\^ace à une approximation forte de la martingale $M$ 
par une trajectoire Brownienne réalisée en exploitant la méthode de troncature. Alors que les résultats  de \cite{f&f&a} 
ont été obtenus en reprenant la technique de la fonction caractéristique utilisée par \cite{touati1} pour démontrer le 
théorème de la limite centrale généralisé pour les martingales. 

Le but de cet article consiste d'une part à  généraliser  le théorème de la limite centrale presque-sûre
 aux  martingales quasi-continues à gauches et d'autre part à établir des théorèmes limites précisant  les vitesses 
de convergences (en loi et au sens presque-sûr) de la loi forte quadratique (LFQ) associée à ce théorème de 
la limite centrale presque-sûre pour les martingales 
 quasi-continues à gauches. 
L'exemple suivant met en évidence l'application des différents théorèmes obtenus et leur usage en statistique
\subsection{Estimation de la variance d'un P.A.I.S.}
Soit $(S_t)_{t\geq 0}$ un processus à accroissements indépendants et stationnaires (P.A.I.S.) dont la mesure de
Lévy des sauts $\nu$ vérifie:
\begin{equation}\label{pais}
\nu(dt,dx)=dt\,F(dx),\;\;\mbox{avec}\;\;\int|x|^{2p}F(dx)<\infty\;\;\mbox{pour un}\;\;  p>1,
\end{equation}
o\`u $F$ est une mesure positive sur $\mathbb R$. 
On note:
$$m=\mathbb E\,S_1 ,\;\;\;\;\sigma^2=\mathbb E\, S_1^2-m^2.$$
La  loi forte quadratique (voir Théorème \ref{LFQ}) nous permet de définir un estimateur fortement consistant de 
$\sigma^2$. En effet on a le résultat suivant    
\begin{equation*}
{\hat \sigma_t}^2:=(\log(1+t))^{-1}\int_0^t\frac{(S_r-mr)^2}{(1+r)^2}\,dr\underset{t\rightarrow \infty}{\longrightarrow}
\sigma^2\;\;\;p.s..
\end{equation*}
Si de plus, pour un $\rho>1/2$  on a 
$$ (1+t)^{-1}\sum_{r\leq t}(\Delta S_r)^2-\int_{\mathbb R}|x|^2F(dx)\leq c^{te}\,[\log(1+t)]^{-\rho}\;\;\mbox{p.s.,}$$ 
alors  le théorème de la limite centrale associé à la loi forte quadratique (voir Théorème \ref{TLCLFQ})
nous permet d'établir le résultat suivant
\begin{equation*}
\sqrt{\log(1+t)}({\hat \sigma}^2-\sigma^2)\Rightarrow \mathfrak N(0,4\sigma^4).
\end{equation*}
Ces résultats seront étendus à des P.A.I.S. pondérés. (voir la partie 4 du papier). 
Les principaux résultats sont énoncés au  paragraphe suivant et  démontrés au paragraphe 4. Au paragraphe 3, on regroupe les outils 
techniques utilisés dans les preuves. Ces outils sont établis au paragraphe 5.    
\section{Préliminaires}
On note  $\|.\|$ la norme Euclidienne sur $\mathbb R^d$. Pour une matrice 
réelle carrée $A$: $A^*$, $\tr(A)$, et $\det(A)$ désignent respectivement la 
matrice transposée, la trace et le déterminant de $A$. La norme de $A$ est 
définie par: $\|A\|^2=\tr(A^*A)$. On considère une martingale quasi-continue 
à gauche
$ M=(M_t)_{t\geq 0}$ $d$-dimensionnelles, localement de carré intégrable, 
définies sur un espace de probabilité filtré $(\Omega,\mathcal F,(\mathcal F_t)_{t\geq 0},\mathbb P)$.
On considère de même 
un processus déterministe $V=(V_t)_{t\geq 0}$ à valeurs dans l'ensemble des matrices inversibles.
Dans la suite on rappelle un théorème fondamental de \cite{touati2} qui nous sera
utile dans les preuves de nos principaux résultats.\\

 Pour $u\in \mathbb R^d$ on définit
\begin{multline*}
\Phi_t(u):=\exp\Bigl(-\frac{1}{2}u^*\langle M^c\rangle_t u\\
+\int_0^t\int_{\mathbb R^d}^{}\bigl(\exp(i\langle u,x)-1-i\langle u,
x\rangle\bigr)\nu^{\s M}(ds,dx)\Big)
\end{multline*}
o\`u $M^c$, $\nu^{\s M}$ sont respectivement la partie martingale continue et la mesure de 
Lévy des sauts de $M$.
\begin{theoreme}[Théorème Limte Centrale Généralisé pour les Martingales]\label{AT}
Soit $M=(M_t)_{t\geq 0}$ une martingale locale, $d$-dimensionnelle, nulle en $0$ 
et quasi-continue à 
gauche. Soit 
$V=(V_t)_{t\geq 0}$ une famille déterministe de matrices inversibles. Si le couple $(M,V)$ vérifie 
l'hypothèse:
$$ (\mathcal H)\;\left\{
\begin{array}{rl}
 \Phi_t((V_t^*)^{-1}u)\rightarrow \Phi_\infty(\eta,u)\;\;p.s.\\\\
 \Phi_{\infty}(\eta,u)\;\;\;\mbox{non nulle}\;\;p.s.\;\;\;\;\;\;\;\;
\end{array}
\right.$$
(o\`u $\eta$ désigne une v.a. sur $(\Omega,\mathcal F,\mathbb P)$, 
éventuellement dégénérée et à valeurs dans un espace vectoriel de dimension finie
$\mathfrak X$)
alors on a
$$Z_t:=V_t^{-1}M_t\Rightarrow Z_{\infty}:=\Sigma(\eta)$$ 
de manière stable o\`u $(\Sigma(x),x\in\mathfrak X)$ est un processus de loi $\mathcal Q$ et 
indépendant de la v.a $\eta$.
\end{theoreme}
Notons  que pour $(x,u)\in\mathfrak X\times\mathbb R^d$:
$$\Phi_{\infty}(x,u)=\int_{\mathbb R^d}^{}\exp\bigl(i\langle u,\xi\rangle\bigr)\pi(x,d\xi)$$
désigne la transformée de Fourier des lois marginales unidimensionnelles 
$\bigl( \pi(x,.);x\in\mathfrak X\bigr)$ d'une loi de probabilité $\mathcal Q$ 
sur l'espace $\mathfrak C(\mathfrak X,\mathbb R^d)$ des fonctions continues de 
$\mathfrak X$ dans $\mathbb R^d.$





\section{Énoncé des principaux résultats}
Dans la suite, on donne quelques propriétés aux quelles doit obéir la normalisation matricielle $(V_t)$.
On dit que la famille $(V_t)$ vérifie la condition $(\mathcal C)$ si les trois 
propriétés $\{(\mathcal C1),\;(\mathcal C2),\;(\mathcal C3)\}$ ont lieu:\\

{\bf $\bullet\,(\mathcal C1)$} $t\,\mapsto\,V_t$ est de classe $\mathscr C^1$;\\

{\bf $\bullet\,(\mathcal C2)$} il existe  $s_0\geq 0$ tel que pour tout  $t\geq s\geq s_0$ on a $V_sV_s^*\leq V_tV_t^*$ 
(au sens des matrices réelles symétriques positives);\\

{\bf $\bullet\,(\mathcal C3)$} il existe une fonction $a=(a_t)$ continue, 
décroissante  vers 0 à l'infini, telle que :$$A_t:=\int_0^ta_sds\uparrow\infty\;\; \mbox{pour}\;\; t\uparrow \infty $$
et une matrice $U$ vérifiant:$$a_t^{-1}V_t^{-1}\frac{dV_t}{dt}-U=\Delta_t,\;\;\;\mbox{avec}\;
\lim_{t\rightarrow\infty}\Delta_t =0$$ 
et telle que  la matrice symétrique $S:=U+U^*$ soit  définie positive.

\subsection{Théorème de la  limite centrale presque-sûre généralisé.}
\begin{theoreme}\label{TLCPSQ}
Soit $M=(M_t)_{t\geq 0}$ une martingale locale, $d$-dimensionnelle, nulle en $0$ 
et quasi-continue à 
gauche. Soit 
$V=(V_t)_{t\geq 0}$ une famille déterministe de matrices inversibles  satisfaisant
 aux conditions $(C)$. 
Si le couple $(M,V)$ vérifie l' hypothèse { $(\mathcal H)$} et l'hypothèse 
$$(\mathcal H1):\;\;V_t^{-1}\langle M\rangle_t(V_t^*)^{-1}\rightarrow C\;\;\;\mbox{p.s.,}\;\;$$ 
({o\`u $C$ 
est une matrice aléatoire ou non}) 
alors les mesures $\left( \mu _{R}\right) ${\
al\'{e}atoires d\'{e}finies par:}
$$\mu _{R}\,=\,\left( \log \text{\thinspace }\left( \det V_{R}^{2}\right)
\right) ^{-1}\int\nolimits_{0}^{R}\,\delta _{Z_{r}}d\log 
\left( \det V_{r}^{2}\right),\;\;\mbox{o\`u }\;\;Z_{r}=V_{r}^{-1}M_{r} $$
{v\'{e}rifient la
version g\'{e}n\'{e}ralis\'{e}e suivante du TLCPS:}

$$\mbox{ (TLCPSG)}:\;\;\;\;\;\;\;
\mu _{R}\Longrightarrow \mu _{\infty }\,\,\,\, \mbox{p.s..}$$
\end{theoreme}
\begin{remarque}
Notons que sous l'hypothèse $(\mathcal H1)$ et l'hypothèse 

$$(\mathcal H'):\;\forall\; \delta> 0,\;\;\;\int_{\mathbb R^d}^{}\int_0^t\|V_t^{-1}x\|^2 
\mathfrak  1_{\{\|V_t^{-1}x\|>\delta\}}\nu^{\s M}(ds,dx)\rightarrow 0.$$
 l'hypothèse $(\mathcal H)$ a lieu avec 
$$\eta=C^{1/2}\;\;\;\mbox{ et }\;\;\;\Phi_{\infty}(x,u)=\exp(-\frac{1}{2}u^*xx^*u).$$
L'hypothèse $(\mathcal H')$ est plus connue sous le nom de condition de Lindberg.
\end{remarque}
\subsection{Lois fortes quadratiques  associées au TLCPS}
Le théorème suivant donne une loi forte des grands nombres  avec une normalisation matricielle:
\begin{theoreme}\label{LFQ}
Soit  $M=(M_t)_{t\geq 0}$ une martingale locale, $d$-dimensionnelle, quasi-continue à 
gauche et nulle en $0$. On suppose que pour une famille de matrices inversibles $V=(V_t)_{t\geq 0}$ vérifiant 
la condition $(\mathcal C)$. Si le couple $(M,V)$ satisfait aux hypothèses: $(\mathcal H)$, 
$(\mathcal H1)$, 

$$(\mathcal H2):\;V_t^{-1}[M]_t(V_t^*)^{-1}\rightarrow C\;\;\;\mbox{p.s..}$$ et 
$$(\mathcal H3):\;C=\int_{}^{} xx^*d\mu_{\infty}(x).$$
(où $\mu_{\infty}=\mu_{\infty}(\omega,.)$ désigne la probabilité de transition (éventuellement non aléatoire)
loi de la v.a $\Sigma\bigl(\eta(\omega)\bigr)$ (voir Théorème \ref{AT})).
 Alors on a les résultats suivants:
 \begin{equation*}
{\mbox{\ (LFQ):}}\;\;\;\; \bigl(\log{(\det V_R^2)}\bigr)^{-1}\int_{0}^{R}V_s^{-1}M_{s_-}M^{*}_{s_-}{V_s^*}^{-1}
d\bigl(\log (\det V_s^2\bigr)\rightarrow C\;\;\;\; p.s.,
\end{equation*} 

\begin{equation*}
\mbox{ (LL):}\;\;\;\;\|V_r^{-1}M_r\|=o\bigl(\sqrt{\log(\det V^2_r)}\bigr)\;\;\;p.s..
\end{equation*}
\end{theoreme}
\begin{remark}
Notons que l'hypothèse $(\mathcal H3)$ est automatiquement vérifiée sous les hypothèses 
$(\mathcal H')$ et $(\mathcal H1)$.
\end{remark}
 \subsection{Vitesses de convergence de  la LFQ (cas d'une normalisation matricielle)}
Dans la suite on donne  un TLC pour la  LFQ établie ci-dessus. 
\begin{theoreme}\label{TLCLFQ}
{Soit M=}$\left( M_{t},\,t\geq 0\right) \,\,$ une martingale locale, $d$-dimensionnelle,
quasi-continue \`{a} gauche, nulle en $0$. On suppose que pour une famille de matrices inversibles $V=(V_t)_{t\geq 0}$ 
vérifiant la condition $(\mathcal C)$ et que le couple $(M,V)$ satisfait
aux hypothèses: $(\mathcal H)$, $(\mathcal H1)$, $(\mathcal H2)$ et $(\mathcal H3)$. Supposons de plus, 
 que  la condition { $(\mathcal C3$)} est vérifiée avec $\Delta_t=O(A_t^{-3/2})$, $(t\rightarrow \infty).$
 Alors si  $R$ désigne la matrice symétrique, positive solution de l'équation de Lyapounov:
$$I=RU+U^*R,$$ on obtient: 
\begin{multline}
\bigl(\log(\det V_t^2)\bigr)^{-1/2}\int_0^t\tr\Bigl[V_s^{-1}\bigl(M_{s_-}M_{s_-}^*-[M]_s\bigr)(V_s^{*})^{-1}\Bigr]
d\bigl(\log(\det V_s^2)\bigr)\hspace{0.25cm}\\\Rightarrow   2\sqrt{\tr(S)\tr\bigl(\tilde CRCR\bigr)}G,
\end{multline}
 o\`u  $\tilde C:= UC+CU^*$, $S=U+U^*$ (U étant la matrice définie dans la condition ($\mathcal C$)) et
 $G$  une gaussienne centrée réduite.
Si de plus pour un $\rho>1/2$ on a que:   
$$\log\bigl(\det(V^2_t)\bigr)^{^{\rho}}\Bigl|\tr\bigl\{V_t^{-1}\bigl([M]_t\bigr)(V_t^{*})^{-1}-C\bigr\}\Bigr|=O(1)\;
\;p.s.,\;\;\;(t\rightarrow \infty)$$
alors 
\begin{multline}
\bigl(\log(\det V_t^2)\bigr)^{-1/2}\int_0^t\tr\Bigl[V_s^{-1}\bigl(M_{s_-}M_{s_-}^*\bigr)(V_s^{*})^{-1}-C\Bigr]
d\bigl(\log(\det V_s^2)\bigr)\\\Rightarrow   2\sqrt{\tr(S)\tr\bigl(\tilde CRCR\bigr)}G.
\end{multline}
\end{theoreme}
Dans ce cas, on donne une loi du logarithme itéré logarithmique associée à la (LFQ) qu'on notera : (LILL) 
\begin{theoreme}\label{LLIL}
{Soit M=}$\left( M_{t},\,t\geq 0\right) \,\,$ une martingale locale, $d$-dimensionnelle,
quasi-continue \`{a} gauche, nulle en $0$.  
On suppose que pour une famille de matrices $V=(V_t)_{t\geq 0}$ vérifiant la condition $(\mathcal C)$, le couple $(M,V)$ satisfait
aux hypothèses: $(\mathcal H)$, $(\mathcal H1)$, $(\mathcal H2)$ et $(\mathcal H3)$. Supposons de plus  que:
$$\mathbb E\bl[\sup_t M_t(\Delta M_t)^*\br]< \infty.$$

On considère $R$ la matrice symétrique, positive solution de l'équation de Lyapounov:
$$I=RU+U^*R.$$ 
Si  pour un $\rho>1/2$ on a que:   
$$\log\bigl(\det(V^2_t)\bigr)^{^{\rho}}\Bigl|\tr\bigl\{V_t^{-1}[M]_t(V_t^{*})^{-1}-C\bigr\}\Bigr|=O(1)\;
\;p.s.,\;\;\;(t\rightarrow \infty)$$
alors 
\begin{multline*}
{{\underset{t\rightarrow \infty}{\overline{\lim}}}}\;\frac{1}{h\bigl(\log(\det V_t^2)\bigr)}\int_0^t \tr\Bigl
\{ V_s^{-1}\bigl(M_{s_-}M_{s_-}^*\bigr)(V_s^{*})^{-1}-C  \Bigr\}\,d\bigl(\log(\det V_s^2)\bigr)\\ \leq \sqrt{\tr(S)
\tr\bigl(\tilde CRCR\bigr)}\;\;\;\mbox{p.s.,}
\end{multline*}
o\`u  $\tilde C:= UC+CU^*$, $S=U+U^*$ (U étant la matrice définie dans la condition ($\mathcal C$))
et $h(u)=\sqrt{2u\log\log\, u}\;\;\mbox{pour\;} u\geq e$.
\end{theoreme}
\subsection{Vitesses de convergence de  la LFQ (cas d'une normalisation scalaire)}
On dira que $V_t$ est une normalisation scalaire vérifiant  la condition ($\mathcal C$) si elle 
est de la forme 

$$V_t=v_tI_d$$ o\'u $v_t$ est une fonction scalaire de classe $\mathcal C^1$ satisfaisant au deux conditions suivantes\\
$\bullet$ il existe un $s_0\geq 0$ tel que  pour tout $t\geq s\geq s_0$ on a  $v_s^2\leq v_t^2$\\
$\bullet$ il existe une fonction $a=(a_t)$ continue, décroissante vers $0$ à l'infini, telle que:
 
$$A_t=\int_0^ta_s\,ds\uparrow \infty \;\;\mbox{pour}\;\;t\uparrow \infty$$ et un scalaire positif $\eta$
tel que  
\begin{equation}\label{speed}
a_t^{-1}v_t^{-1}v'_t-\eta =\delta_t,\;\;\mbox{avec}\;\; \lim_{t\rightarrow \infty}\delta_t=0.
\end{equation}
Ainsi le théorème limite centrale asocié à la (LFQ) est donné par le résultat suivant 
\begin{theoreme}\label{tlclfqs}
{Soit M=}$\left( M_{t},\,t\geq 0\right) \,\,$ une martingale locale, $d$-dimensionnelle,
quasi-continue \`{a} gauche, nulle en $0$. Soit $V=(V_t)_{t\geq 0}$ une normalisation scalaire vérifiant la condition $(\mathcal C)$ et tel que 
 le couple $(M,V)$ satisfait
aux hypothèses: $(\mathcal H)$, $(\mathcal H1)$, $(\mathcal H2)$ et $(\mathcal H3)$. Supposons de plus  
 que  la relation { (\ref{speed})} est vérifiée avec $\delta_t=O(A_t^{-3/2})$, $(t\rightarrow \infty).$
 Alors il vient que
\begin{equation}
\bigl(\log(v_t^2)\bigr)^{-1/2}\int_0^t v_s^{-2}\bigl[M_{s_-}M_{s_-}^*-[M]_s\bigr]
d\bigl(\log(v_s^2)\bigr)\hspace{0.25cm}\Rightarrow   (2\eta C)\,G,
\end{equation}
 o\`u  $G$ est une gaussienne centrée réduite.
Si de plus on suppose que  pour un $\rho>1/2$ on a    
$$\log\bigl(v^{2\rho}_t\bigr)\Bigl|v_t^{-2}[M]_t-C\Bigr|=O(1)\;
\;p.s.,\;\;\;(t\rightarrow \infty)$$
alors 
\begin{equation}
\bigl(\log(v_t^2)\bigr)^{-1/2}\int_0^t \bigl[v_s^{-2}M_{s_-}M_{s_-}^*-C\bigr]
d\bigl(\log( v_s^2)\bigr)\\\Rightarrow   (2\eta C)\,G.
\end{equation}
\end{theoreme}
Dans ce cadre la loi du logarithme itérée est donnée par le théorème suivant
\begin{theoreme}\label{llils}
{Soit M=}$\left( M_{t},\,t\geq 0\right) \,\,$ une martingale locale, $d$-dimensionnelle,
quasi-continue \`{a} gauche, nulle en $0$.  
On suppose que pour une normalisaion scalaire $V=(V_t)_{t\geq 0}$ vérifiant la condition $(\mathcal C)$,
 le couple $(M,V)$ satisfait
aux hypothèses: $(\mathcal H)$, $(\mathcal H1)$, $(\mathcal H2)$ et $(\mathcal H3)$. Supposons de plus  que:
$$\mathbb E\bl[\sup_t M_t(\Delta M_t)^*\br]< \infty.$$

Si  pour un $\rho>1/2$ on a que:   
$$\log\bigl(v^{2\rho}_t\bigr)\Bigl|v_t^{-2}[M]_t-C\bigr\}\Bigr|=O(1)\;
\;p.s.,\;\;\;(t\rightarrow \infty)$$
alors 
\begin{equation*}
{{\underset{t\rightarrow \infty}{\overline{\lim}}}}\;\frac{1}{h\bigl(\log( v_t^2)\bigr)}\int_0^t 
\bigl[v_s^{-2}M_{s_-}M_{s_-}^*-C\bigr]\,d\bigl(\log( v_s^2)\bigr)\\ \leq 2\eta C\;\;\;\mbox{p.s.,}
\end{equation*}
o\`u $h(u)=\sqrt{2u\log\log\, u}\;\;\mbox{pour\;} u\geq e$.
\end{theoreme}



\section{Démonstration des principaux résultats}
Au début de ce paragraphe, on donne une propriété simple nous permettant de simplifer
les preuves des principaux  résultats. En effet, 
on rappelle que la differentielle du determinant d'une matrice inversible $X$ 
est donnée par
\begin{equation}\label{derive}
d\det(X)=\det(X)\tr(X^{-1}dX)
\end{equation}
On en déduit alors que 
\begin{equation}\label{EQ1}
\int_0^t2\tr\bigl[V_s^{-1}\frac{dV_s}{ds}\bigr]ds=\log (\det V_{t})^{2},
\end{equation}
Compte tenu des conditions  $(\mathcal C)$ on voit que 
\begin{equation}\label{equiv}
\frac{\log\bigl(\det(V^2_t)\bigr)}{A_t\tr(S)}\rightarrow 1\;\;\;\mbox{p.s.}
\end{equation}
avec $S=U+U^*$ la matrice introduite dans $(\mathcal C_3)$. Ainsi, cette propriété permettera de remplacer certaines moyennes logarithmiques par des moyennes pondérées par la fonction $a$.

\subsection{Preuve du Théorème \ref{TLCPSQ}}
Pour démontrer le Théorème \ref{TLCPSQ} on va étudier la fonction caractéristique associée aux mesures $(\mu_R)$
donnée par 
$$\psi_R(u)=\,\left( \log \text{\thinspace }\left( \det V_{R}^{2}\right)
\right) ^{-1}\int\nolimits_{0}^{R}\,\exp\{i\langle u, {Z_{r}}\rangle \}d\log 
\left( \det V_{r}^{2}\right)$$
En vue de simplifier la preuve on démontre tout  d'abord le lemme suivant
\begin{lemme}
Sous les hypothèses du Théorème \ref{TLCPSQ}
$$\lambda_R:=\psi_R(u)-A_{R}^{-1}\int
\nolimits_{0}^{R}\exp \left\{ i\left\langle u,Z_{r}\right\rangle \right\}
dA_{r}\rightarrow 0\;\;p.s.\;\;(R\rightarrow 0).$$
\end{lemme}
\begin{proof}[Preuve]
En décomposant l'expression de $\lambda_R$ comme suit 
$$\lambda_R=\lambda_R^1+\lambda_R^2,$$
avec

$$\lambda^1_R:=\left( \log \text{\thinspace }\left( \det V_{R}^{2}\right)
\right) ^{-1}\int\nolimits_{0}^{R}\,\exp\{i\langle u, {Z_{r}}\rangle \}d\bigl(\log 
\left( \det V_{r}^{2}\right)- \tr(S)A_r\bigr)$$
et
$$
\lambda_R^2:=\bigl(\left( \tr(S)\log \text{\thinspace }\left( \det V_{R}^{2}\right)
\right) ^{-1}-A_R^{-1}\bigr)\int
\nolimits_{0}^{R}\exp \left\{ i\left\langle u,Z_{r}\right\rangle \right\}
dA_{r},$$
en remarquant que 
$$|\lambda_R^1|\leq\bigl| \log \text{\thinspace }( \det V_{R}^{2})\bigr|^{-1}
\bigl|\log 
\left( \det V_{R}^{2}\right)- \tr(S)A_R\bigr|$$
on déduit par la relation (\ref{equiv}) que $\lambda_R^1\rightarrow 0,\,(R\rightarrow 0)$.
De la même façon on voit que 
$$|\lambda_R^2|\leq\Bigl|\frac{A_R\tr(S)}{\log 
\left( \det V_{R}^{2}\right)}-1\Bigr|\rightarrow 0.$$
Ce qui termine la preuve du lemme.
\end{proof}
Compte tenu du lemme précédent on conclut que pour démontrer la propriété (TLCPS) il nous  suffit de prouver que 
\begin{equation}\label{aim} 
A_{R}^{-1}\int
\nolimits_{0}^{R}\exp \left\{ i\left\langle u,Z_{r}\right\rangle \right\}
dA_{r}\rightarrow \Phi _{\infty }(\eta ,u)
\end{equation}
avec $\Phi _{\infty }(\eta ,u)$ est la fonction caractéristique associée à la mesure  limite $\mu_{\infty}$.
Ainsi, en vue de démonter cette dernière relation, on va  expliciter  l'expression  de la variable aléatoire 
complexe $\exp \left\{ i\left\langle u,Z_{r}\right\rangle \right\}$. 
Pour ce fait, On rappllera  quelques résultats utiles dans  la suite. On note
$ \Phi_t(u):=\exp\{ B_t(u)\}$
avec 
\begin{multline*}
B_{t}\left( u\right) :=-\dfrac{1}{2}\,\,
u^*\left\langle M^{c}\right\rangle _{t}\,\,u\,\\
+\int\nolimits_{0}^{t}\int\nolimits_{\Bbb{R}^{d}}\Bigl( \exp \left\{
i\left\langle u,x\right\rangle \right\} -1-i\left\langle
u,x\right\rangle \Br) \nu^{\s M} \left( ds,dx\right) .
\end{multline*}
Soit $(L_t(u))_{t\geq 0}$  le processus défini par 
$$ L_t(u):= \left[ \Phi _{t}(u)\right] ^{-1}\exp i\left\langle
u,M_{t}\right\rangle, $$ alors on a le résultat suivant
\begin{lemme}
Le processus $(L_t(u))_{t\geq 0}$ est une martingale locale complexe. De plus on a
\begin{equation}\label{mart}
\left| L_{t}\left( u\right)
\right| \leq \exp \left\{ \,\frac{1}{2}u^*\,\left\langle M
\right\rangle_{t}\, u \right\}.
\end{equation}
\end{lemme}
\begin{proof}[Preuve.]
Comme 
$\left( B_{t}\left( u\right)\right)_{t\geq
0}$ est un processus continu on en déduit que \\ $(L_t(u))_{t\geq 0}$ est une martingale 
locale complexe (cf. \cite{jacod&shirayev}).
D'autre part on voit que son  module vaut
\begin{multline}\label{majoration}
\left| L_{t}\left( u\right) \right|
=\exp \left\{ \dfrac{1}{2}\,\,u^*\,\,\left\langle
M^{c}\right\rangle _{t}\,\,u\right\}\\
\times\exp \left\{\int\nolimits_{0}^{t}\int\nolimits_{\Bbb{R}^{d}}\left( 1-\cos
\left\langle u,x\right\rangle \right) \nu^{\s M} \left( ds,dx\right)
\right\}.
\end{multline}
Ainsi la majoration (\ref{mart}) découle directement du fait que  
$$1-\cos x\leq x^2/2,\;\;\forall \,x\in \mathbb R.$$
\end{proof}
Par conséquent, démontrer la relation (\ref{aim}) revient à prouver que 
\begin{equation}\label{cible}
A_{R}^{-1}\int
\nolimits_{0}^{R}L_r\bigl((V_r^*)^{-1}u\bigr)\Phi_r\bigl((V_r^*)^{-1}u\bigr)
dA_{r}\rightarrow \Phi _{\infty }(\eta ,u).
\end{equation}
Ainsi, afin d'exploiter le lemme précédent on introduit les temps d'arrêts suivants.
Pour $u$\ \ fix\'{e} dans $\Bbb{R}^d$, soit  $b>0$ un point de continuit\'{e}
de la v.a. $\tr(C)$ et $c>0$ un point de continuité de la v.a. 
$|\Phi_{\infty}(\eta,u)|^{-1}$.
Considérant les événements
$$E_{r}^{b}=\left\{ \tr(C_{r}) >b\right\}\;\;\mbox{ et }\;\; 
E_{r}^{u,c}=\left\{ |\Phi_r(u)|^{-1} >c\right\},$$ o\`u 
$$C_r:=V_{r}^{-1}\langle
M\rangle_{r}\,\,(V_{r}^*)^{-1},$$ on définit  le temps d'arrêt:
 
$$T_r :=T_{r}^{b,c}\left( u\right)=T_r^b\wedge T_r^c(u),$$
avec
\begin{align*}
 T_r^b:=\left\{ 
\begin{array}{l}
\inf\left\{ t\leq r\;\;/\;\;\tr\bl(V_{r}^{-1}\langle
M\rangle_{t}\,\,(V_{r}^*)^{-1}\br)>b\right\}\;\;\;\;\; \mbox{si }\;\;\;\;
E_{r}^{b}\mbox{ est r\'{e}alis\'{e},}\\\\ 
r\hspace{1cm}\mbox{sinon}
\end{array}
\right.
\end{align*}
et 
\begin{align*}
 T_r^c(u):=\left\{ 
\begin{array}{l}
\inf\left\{ t\leq r\;\;/\;\;|\Phi_{t}\bigl((V_r^*)^{-1}u\bigr)|^{-1}>c\right\}\;\;\;\;\; \mbox{si }\;\;\;\;
E_{r}^{u,c}\mbox{ est r\'{e}alis\'{e},}\\\\ 
r\hspace{1cm}\mbox{sinon}
\end{array}
\right.
\end{align*}
Notons que d'apr\`{e}s l'in\'{e}galit\'{e} (\ref{mart}), $\Bigl( L_{t\wedge T_{r}}\bigl((V_r^*)^{-1}
u\bigr)\Bigr)_{t\geq 0} $ est une martingale locale complexe dont le module
est major\'{e} par $\exp(b\left\| u\right\|^{2}/2).$ C'est
donc une martingale d'esp\'{e}rance 1. D'o\`{u} la propri\'{e}t\'{e}: $$
\Bbb{E}\; L_{r\wedge T_{r}}\bigl((V_r^*)^{-1}
u\bigr) =1.$$
Il vient alors que
\begin{multline}
A_{R}^{-1}\int
\nolimits_{0}^{R}L_r\bigl((V_r^*)^{-1}u\bigr)\Phi_r\bigl((V_r^*)^{-1}u\bigr)
dA_{r} - \Phi _{\infty }(\eta ,u)=\\ A_{R}^{-1}\int
\nolimits_{0}^{R}\Bigl[L_{r\wedge T_r}\bigl((V_r^*)^{-1}u\bigr) -1\Bigr] \Phi _{\infty }(\eta ,u)dA_r+ 
\Delta _{R}\left( b,c,u\right)
\\+\delta'_R(b,c,u) + \delta''_R(b,c,u)
\end{multline}
avec 
\begin{multline}
\Delta _{R}\left( b,c,u\right):= A_{R}^{-1}\int_{0}^{R}\exp \left\{
i\left\langle u,V^{-1}_rM_r\right\rangle \right\}
dA_{r}\\-A_{R}^{-1}\int_{0}^{R}\exp \left\{ i\left\langle
u,V_r^{-1}M_{r\wedge T_{r}}\right\rangle \right\} dA_{r},
\end{multline}
$$\delta'_R(b,c,u):=A_{R}^{-1}\int\nolimits_{0}^{R}L_r\bigl((V_r^*)^{-1}u\bigr) \bl[\Phi_r\bigl((V_r^*)^{-1}u\bigr)
-\Phi _{\infty }(\eta ,u)\br] dA_{r},$$
et
$$\delta''_R(b,c,u):=A_{R}^{-1}\int\nolimits_{0}^{R}L_r\bigl((V_r^*)^{-1}u\bigr) \bl[ \Phi_{r\wedge T_r}\bigl((V_r^*)^{-1}u\bigr)
-\Phi_r\bigl((V_r^*)^{-1}u\bigr)\br] dA_{r}.$$
Par conséquent la relation (\ref{cible}) est immédiate dés que les deux propriétés suivantes sont vérifiées 
\begin{equation}\label{delta0}
\Delta _{R}\left( b,c,u\right)
\\+\delta'_R(b,c,u) + \delta''_R(b,c,u)\rightarrow 0
\end{equation}
et 
\begin{equation}\label{marta}
 A_{R}^{-1}\int
\nolimits_{0}^{R}\Bigl[L_{r\wedge T_r}\bigl((V_r^*)^{-1}u\bigr) -1\Bigr]dA_r\rightarrow 0\,\,p.s..
\end{equation}
\subsubsection{Vérification de la propriété (\ref{delta0})}
Comme  $L_r\bigl((V_r^*)^{-1}u\bigr)\leq c$ on en déduit que 
$$|\delta'_R(b,c,u)|\leq c\, A_{R}^{-1}\int_{0}^{R}\bl|\Phi_r\bigl((V_r^*)^{-1}u\bigr)-
\Phi _{\infty }(\eta ,u)\br| dA_{r}.$$
Ainsi vu l'hypothèse $(\mathcal H)$ il vient que  $$\underset{R\rightarrow \infty 
}{\overline{\lim }}\left| \delta' _{R}\left( b,c,u\right) \right|
\longrightarrow 0 \mbox{ p.s..}$$
Par ailleurs, on voit que  
\begin{equation*}
\Delta_R(b,c,u)\vee \delta''_R(b,c,u)\leq 2cA^{-1}_R\int_0^R\mathbf 1_{\{T_r<r\}}dA_r.
\end{equation*}
Or on sait que d'une part
\begin{align*}
\mathbf 1_{\displaystyle\{T_r<r\}}&\leq \mathbf 1_{\displaystyle \{T^b_r<r\}}+\mathbf 1_{\displaystyle\{T^{c}_r(u)<r\}}\\
                     & \leq \mathbf 1_{\displaystyle E^b_r}+\mathbf 1_{\displaystyle E^{u,c}_r}
\end{align*}
et que d'autre part  $\mathbb P(\tr(C)=b)=\mathbb P(|\Phi_{\infty}(\eta,u)|^{-1}=c)=0$. Par conséquent, à l'aide des hypothèses 
$(\mathcal H1)$ et $(\mathcal H)$ il vient que 
$$ \underset{R\rightarrow \infty}{\lim \sup}\,
\Delta_R(b,c,u)\vee \delta''_R(b,c,u)\leq 2 c\Bl(\mathbf 1_{\displaystyle \{\tr(C)>b\}}+\mathbf 1_{\displaystyle\{|\Phi_{\infty}(\eta,u)|^{-1}>c\}}\Br).$$
Ainsi en faisant tendre $b$ et $c$ de manière séquentielle et de sorte qu' on ait toujours 
$\mathbb P(\tr(C)=b)=\mathbb P(|\Phi_{\infty}(\eta,u)|^{-1}=c)=0$, on obtient que  
$$\underset{R\rightarrow \infty 
}{\overline{\lim }}\left| \delta _{R}\left( b,c,u\right) \right|
\longrightarrow 0 \mbox{ p.s..}$$
En vue de simplifier les notations on pose
$$ \tilde L_r(u):=L_{r\wedge T_r}\bigl((V_r^*)^{-1}u\bigr).$$
Le reste de la preuve du théorème, consiste à \'{e}tablir la convergence p.s. des moyennes $A_{R}^{-1}\int%
\nolimits_{0}^{R}\tilde L_r(u)dA_{r}$ vers 1. On se propose alors de  montrer d'abord que
cette convergence \`{a} lieu en moyenne quadratique. D'o\`{u} l'\'{e}tape
cruciale suivante consacr\'{e}e \`{a} l'estimation de la covariance du
couple $ \bigl(\tilde L_r(u),\tilde L_{\rho}(u)\bigr) .$
\subsubsection{Estimation de la covariance du couple $\left(
\tilde L_{r}\left( u\right) , \tilde L_{\rho }\left( u\right) \right) .$}

Pour tous\textit{\ }$u\in \Bbb{R}^d,$ $\left( \rho ,r\right) \in \Bbb{R}
_{+}\times \Bbb{R}_{+}$ \ avec $\rho \leq r,$ notons

$$K_{\rho ,r}\left( u\right) :=\Bbb{E}\left\{ \left( \tilde L_{\rho }\left( u\right)
-1\right) \overline{\left( \tilde L_{r }\left( u\right) -1\right) }\right\}.$$

Comme $\left( L_{r,t\wedge T_{r}}\left(
u\right)\right)_{t\geq0} $ est une martingale  on vérifie aisément  que:
\begin{align*}
K_{\rho ,r}\left( u\right) &=\Bbb{E}\left\{ \tilde L_{\rho }\left( u\right) 
\overline{ \tilde L_{r}\left( u\right) }\right\} -1\\&=\Bbb{E}\left\{ \tilde L_{\rho }\left(
u\right) \overline{\Bbb{E}\left\{ \tilde L_{r}\left( u\right) /\frak{F}_{\rho
\wedge T_{\rho}}\right\} }\right\}-1;
\\&=\Bbb{E}\left\{ \tilde L_{\rho}\left(
u\right) \overline{  L_{\rho \wedge T_{\rho}}\bigl((V_r^*)^{-1}u\bigr)}\right\}-1\\
&=\Bbb{E}\left\{ \tilde L_{\rho}\left(
u\right) \bigl[\overline{ L_{\rho \wedge T_{\rho}}\bigl((V_r^*)^{-1}u\bigr)}-1
\bigr]\right\}
\end{align*}
Ainsi l'inégalité de Cauchy Schwarz donne
\begin{align*}
\left| K_{\rho ,r}\left( u\right) \right| &\leq \left( \Bbb{E} \left|
 \tilde L_{\rho }(u)\right| ^{2} \right) ^{1/2}\left( \Bbb{E} \left|
L_{\rho \wedge T_{\rho}}\bigl((V_r^*)^{-1}u\bigr)-1\right| ^{2} \right) ^{1/2}\\
&\leq  \left( \Bbb{E} \left|
L_{\rho }(u)\right| ^{2} \right)^{1/2}
\left(\Bbb{E} \left| L_{\rho \wedge T_{\rho}}\bigl((V_r^*)^{-1}u\bigr)\right|
^{2} -1  \right) ^{1/2}.
\end{align*}
Or de l'in\'{e}galit\'{e} (\ref{mart}) on  voit que  d'une part  
\begin{align*}
\mathbb E\left|
L_{\rho }(u)\right| ^{2} &\leq \mathbb E\exp \left\{
\,u^*\,V_{\rho}^{-1}\left\langle M\right\rangle _{\rho \wedge
T_{\rho}}\,(V_{\rho}^*)^{-1}\,\,u\right\}\\ &\leq \exp \left\{ b\left\|
u\right\| ^{2}\right\}
\end{align*}
et que d'autre part 
\begin{align*}
\Bbb{E} \left| L_{\rho \wedge T_{\rho}}\bigl((V_r^*)^{-1}u\bigr)\right|
^{2} &\leq \mathbb \exp \left\{
\,u^*\,V_{r}^{-1}\left\langle M\right\rangle _{\rho \wedge
T_{\rho}}\,(V_{r}^*)^{-1}\,\,u\right\}\\ &\leq \exp \left\{ b\left\|
u\right\| ^{2}\left\| V_{r}^{-1}V_{\rho }\right\| ^{2}\right\}.
\end{align*}
Ensuite, en utilisant  l'in\'{e}galit\'{e}: $\forall
t>0,\,\,\,e^{t}-1$\bigskip $\leq te^{t}$ il vient que 
\begin{align*}
\Bbb{E} \left| L_{\rho \wedge T_{\rho}}\bigl((V_r^*)^{-1}u\bigr)\right|
^{2} -1&\leq b\left\|
u\right\| ^{2}\left\| V_{r}^{-1}V_{\rho }\right\| ^{2}\exp \left\{ b\left\|
u\right\| ^{2}\left\| V_{r}^{-1}V_{\rho }\right\| ^{2}\right\}.
\end{align*}
La preuve du  lemme suivant est  explicitée dans la  dernière
section.
\begin{lemme}\label{chaab3}
{ Si la normalisation }$\left(
V_{r}\right) ${ v\'{e}rifie les conditions ($\mathcal{C}$)} 
{alors, pour tout} $\left( r,\rho \right)
\in \mathbb{R}_{+}\times \mathbb{R}_{+}$ {avec }$\rho \leq r$, il existe 
$n_0\in \mathbb N$ tel que

$$\left\| V_{r}^{-1}V_{\rho }\right\| ^{2}\leq d^{n_0}\Bigl(
\frac{\det V_{\rho }}{
\det V_{r}}\Bigr)^{\frac{2}{d}}.$$
\end{lemme}
\vspace{.5cm}
En tenant compte du résultat précédent il vient que 
\begin{equation}\label{encadre}
\left| K_{\rho ,r}\left( u\right) \right|  \leq c^{te}\,\Bigl(\frac{\det V_{\rho }}{
\det V_{r}}\Bigr)^{\frac{2}{d}}
\end{equation}
\noindent pour une constante ind\'{e}pendante de\textit{\ }$\rho $ et de $r$.
\subsubsection{Convergence presque sûre de $
A_{R}^{-1}\int\nolimits_{0}^{R} \tilde L_{r}\left( u\right) dA_{r}$ vers 1.}

 Dans la suite on  v\'{e}rifie d'abord que
\begin{equation}\label{q2}
\Bbb{E}\left\{ \left| \int\nolimits_{0}^{R}\left( \tilde L_{r}(u)-1\right)
d A_r \right| ^{2}\right\} =O\left( A_{R}\right) \left(
R\rightarrow \infty \right),
\end{equation}
En effet, 
$$\Bbb{E}\left\{ \left| \int\nolimits_{0}^{R}\left( \tilde L_{r}(u)-1\right)
d A_r \right| ^{2}\right\} =2 \int_0^R\int_0^r K_{\rho ,r}\,dA_\rho\,dA_r$$
et par l'inégalité (\ref{encadre}), il vient que
\begin{equation}
 \Bbb{E}\left\{ \left| \int\nolimits_{0}^{R}\left( \tilde L_{r}(u)-1\right)
d A_r \right| ^{2}\right\}\leq c^{te}\int_0^R\int_0^r \left|
\frac{\det{V_{\rho}}}{\det{V_r}}\right|^{\frac{2}{d}}\,dA_\rho\,dA_r.
\end{equation}
Or, en utilisant la relation (\ref{EQ1}) on voit que
\begin{align*}
\frac{\det{V_\rho}}{\det{V_r}}&=\exp\Bigl\{-\int_{\rho}^r\tr\bigl[V_s^{-1}
\frac{dV_s}{ds}\bigr]ds\Bigr\}\\
&=\exp\Bigl\{-\int_{\rho}^r\tr\bigl[a_s^{-1}V_s^{-1}
\frac{dV_s}{ds}\bigr]dA_s\Bigr\}
\end{align*}
et donc par la condition {\bf $(\mathcal C3)$} il vient que 
\begin{align*}
\frac{\det{V_\rho}}{\det{V_r}}=\exp\Bigl\{-\tr{[U]}(A_r-A_{\rho})
-\int_{\rho}^r \Delta_s\,dA_s\Bigr\}.
\end{align*}
On en déduit alors qu' il existe $r_0>0$ tel que $\forall r\geq r_0$ on a
\begin{align*}
\int_0^{r}(\frac{\det{V_\rho}}{\det{V_r}})^{\frac{d}{2}}\,dA_{\rho}
&\leq \int_0^r\exp\Bigl\{-\frac{d}{4}\tr{[U]}(A_r-A_{\rho})\Bigr\}\,dA_{\rho}
\\
&\leq \frac{4}{d\tr{U}}\Bigl[1-\exp\Bigl\{-\frac{d}{4}\tr{[U]}A_r\Bigr\}\Bigr]\\
&\leq \frac{4}{d\tr{U}}
\end{align*}
puisque $U$ est une matrice définie positive. D'o\`u le résultat annocé en 
(\ref{q2}). 
Ainsi,  $A_{R}^{-1}\int\nolimits_{0}^{R} \tilde L_{r}\left( u\right) dA_{r}$
tend vers 1 en moyenne quadratique.
Posant $$R_k=\inf\,\{\,r\;/\;\forall\,t>r,\;\;\;A_t>k^2\},$$
il est clair que $A_{R_k}=O(k^2)$ $(k\rightarrow \infty)$. Ainsi il vient que  
\begin{equation*}
\Bbb{E}\left\{ \left| A^{-1}_{R_k}\int\nolimits_{0}^{R_k}\left( \tilde L_{r}(u)-1\right)
dA_r \right| ^{2}\right\} =O\left( k^{-2}\right) \left(
k\rightarrow \infty \right),
\end{equation*}
on en déduit alors que 
$$ A^{-1}_{R_k}\int\nolimits_{0}^{R_k} \tilde L_{r}(u)dA_r\rightarrow 1\;\;\mbox{p.s..}$$
Or pour $R\in[R_k,R_{k+1}[$ on a:
\begin{align*}
\Bl|A^{-1}_{R}&\int\nolimits_{0}^{R} \bl( \tilde L_{r}(u)-1\br)dA_r-A^{-1}_{R_k}
\int\nolimits_{0}^{R_k} \bl( \tilde L_{r}(u)-1\br)dA_r\Br|
\\&\leq \Bl|A^{-1}_{R}\int\nolimits_{0}^{R} \bl( \tilde L_{r}(u)-1\br)dA_r- 
A^{-1}_{R}\int\nolimits_{0}^{R_k} \bl( \tilde L_{r}(u)-1\br)dA_r\Br|\\
&+\Bl|A^{-1}_{R}\int\nolimits_{0}^{R_k} \bl( \tilde L_{r}(u)-1\br)dA_r- 
A^{-1}_{R_k}\int\nolimits_{0}^{R_k} \bl( \tilde L_{r}(u)-1\br)dA_r\Br|
\\&\leq A^{-1}_{R_k}\int\nolimits_{R_k}^{R_{k+1}} \bl(| \tilde L_{r}(u)|+1\br)\,dA_r 
+ |A_R^{-1}-A_{R_k}^{-1}|\int_0^{R_k}\bl(| \tilde L_r(u)|+1\br)\,dA_r
\\&\leq 2(1+c)A_{R_k}^{-1}\bl(A_{R_{k+1}} -A_{R_k}\br)\\
&=O\Bl(\frac{1}{k}\Br)\;(k\rightarrow \infty).
\end{align*}
Puisque $\bl(A_{R_{k+1}} -A_{R_k}\br)=O\Bl({k}\Br)$. Ce qui achève la preuve du  Th\'{e}or\`{e}me
 \ref{TLCPSQ}.
\hfill $\Box$
\subsection{Preuve du Théorème \ref{LFQ}} 
Pour $Z_t:=V_t^{-1}M_t$ et $S=U+U^*$ ($S$ étant
la matrice régulière de la condition $(\mathcal C_3)$), La  
première partie de cette preuve consiste à démontrer la relation suivante

\begin{equation}\label{11}
\displaystyle A_t^{-1}\bigl(\|Z_t\|^2+\int_0^t {Z_{s_-}  }^*SZ_{s_-}  dA_s\bigr)\underset{t
\rightarrow \infty}{\longrightarrow}\tr\bigl(C^{1/2}SC^{1/2}\bigr) \;\;p.s..
\end{equation}
En  effet, en appliquant la formule d'Itô à la  semimartingale $\|Z_t\|^2$ on obtient la
 relation suivante:
\begin{multline}\label{e1}
\|Z_t\|^2=2\int_0^tZ_{s_-}^{*}V_s^{-1}dM_{s_-}  + \tr\bigr(\int_0^t{(V_s^*)}^{-1}V_s^{-1}d[M]_s\bigr)\\
-\int_0^tZ_{s_-}  ^*{V_s}^{-1}d(V_sV_s^*){(V_s^*)}^{-1}Z_{s_-}  ,
\end{multline}
Dans la suite on pose:
$$D_t=\int_0^t Z_{s_-}^*V_s^{-1}d(V_sV_s^*){(V_s^*)}^{-1}Z_{s_-},$$
$$K_t=\int_0^tV_s^{-1}d[M]_s{(V_s^*)}^{-1}\;\;\mbox{ et }\;\;
L_t=\int_0^tZ_{s_-}  ^*V_s^{-1}dM_{s_-}.$$
Avec ces notations, l'égalité (\ref{e1})  s'écrit:
\begin{equation}\label{e3}
\|Z_t\|^2+D_t=2L_t+\tr(K_t)
\end{equation}
et on a les résultats suivants
\begin{lemme}\label{e5}
 \begin{equation}
\frac{K_t}{A_t}\underset{t\rightarrow \infty}{\longrightarrow} CU^*+UC\;\;p.s.. 
\end{equation} 
\end{lemme}
\begin{proof}[Preuve.]
Par la formule d'intégration par parties on voit que
\begin{multline}\label{IPP}
\displaystyle d\bigl(V_s^{-1}[M]_s(V_s^*)^{-1}\bigr)=V_s^{-1}d[M]_s(V_s^*)^{-1}-V_s^{-1}(dV_s)V_s^{-1}[M]_s
(V_s^*)^{-1}\\
                                                    -V_s^{-1}[M]_s(V_s^{*})^{-1}(dV_s)^*(V_s^{*})^{-1},\hspace{2.3cm}
\end{multline}
donc
$$K_t=C'_t+\int_0^tC'_s(V_s^{-1}\frac{dV_s}{ds})^*ds+\int_0^t(V_s^{-1}\frac{dV_s}{ds})C'_sds$$
avec $C'_t:=V_t^{-1}[M]_t{(V_t^*)}^{-1}$.
Par conséquent,
$$K_t=C'_t+\int_0^tC'_s(a_s^{-1}V_s^{-1}\frac{dV_s}{ds})^*dA_s+\int_0^t(a_s^{-1}V_s^{-1}\frac{dV_s}{ds})C'_s\,dA_s.$$
Vu l'hypothèse { $(\mathcal H2)$} et les conditions $(C)$, on déduit  le résultat par lemme de Toeplitz.
\end{proof}
\begin{lemme}\label{e7}
\begin{equation}
D_t\sim\int_0^tZ_s^*SZ_sdA_s\;\;p.s.\;\;(t\rightarrow \infty).
\end{equation}
\end{lemme}
\begin{proof}[Preuve.]
On a 
$$ D_t=\int_0^tZ_{s_-}^*V_s^{-1}d(V_sV_s^*){(V_s^*)}^{-1}Z_{s_-}  =\int_0^tZ_{s}  ^*
\Bigl[\bigl(V_s^{-1}\frac{dV_s}{ds}\bigr)+\bigl(V_s^{-1}\frac{dV_s}{ds}\bigr)^*\Bigr]Z_s\;ds,$$
ainsi compte tenu des conditions $(C)$ et du Lemme de Toeplitz, on déduit aisément le résultat annoncé.
\end{proof}
\begin{lemme}\label{e9}
\begin{equation}
L_t=o(A_t)\;\;p.s..
\end{equation}
\end{lemme}
\begin{proof}[Preuve.]
La variation quadratique prévisible de la martingale locale $(L_t)_{t\geq 0}$ vaut:
$$\langle L\rangle_t=\int_0^t Z_s^*V_s^{-1}d\langle M\rangle_s{(V_s^*)}^{-1}Z_s=\int_0^tZ_{s_-}  ^*d\tilde K_sZ_{s_-}  $$
o\`u $$\tilde K_t=\int_0^t V_s^{-1}d\langle M\rangle_s{(V_s^*)}^{-1},$$ 
est le compensateur prévisible du processus $(K_t)_{t\geq0}$. En utilisant une 
décomposition semblable à celle de $K_t$, on vérifie que 
$$\tilde K_t=C_t+ \int_0^tC_s\bigl(V_s^{-1}\frac{dV_s}{ds}\bigr) ds+
\int_0^t\bigl(V_s^{-1}\frac{dV_s}{ds}\bigr)C_s ds$$
avec $C_t=V_t^{-1}\langle M\rangle_t{(V_t^*)}^{-1}$. Par suite 
\begin{align*}
\langle L\rangle_t &=\int_0^tZ_{s_-}^*dC_sZ_{s_-}  +\int_0^tZ_{s_-}^*\Bigl[C_s\bigl(a_s^{-1}V_s^{-1}\frac{dV_s}{ds}\bigr)^*
+\bigl( a_s^{-1}V_s^{-1}\frac{dV_s}{ds}\bigr)C_s\Bigr]Z_{s_-}  dA_s\\
                   &=\int_0^t \tr (Z_{s_-}Z^*_{s_-}dC_s)\\ &\hspace{2cm}+ \int_0^t \tr\biggl[Z_{s_-}Z_{s_-}^*
\Bigl[C_s\bigl(a_s^{-1}V_s^{-1}\frac{dV_s}{ds}\bigr)^*+\bigl(a_s^{-1}V_s^{-1}\frac{dV_s}{ds}\bigr)C_s\Bigr]\biggr] dA_s\\
&=O\Bigl(\int_0^t \|Z_s\|^2\tr (dC_s)\Bigr) \\&\hspace{2cm}+ O\Biggl(\int_0^t \|Z_s\|^2 \tr\biggl[
\Bigl[C_s\bigl(a_s^{-1}V_s^{-1}\frac{dV_s}{ds}\bigr)^*+\bigl(a_s^{-1}V_s^{-1}\frac{dV_s}{ds}\bigr)C_s\Bigr]\biggr] dA_s\Biggr)
\end{align*}
Vu l'hypothèse { $(\mathcal H1)$},on obtient par le Lemme de Toeplitz: 
$$\langle L\rangle_t= O\Bigl(\int_0^t \|Z_s\|^2\tr (dC_s)\Bigr)  + O(D_t)\;\;p.s..$$
La formule d'intégration par parties et la relation  (\ref{e3}) donnent: 
\begin{align*}
\int_0^t \|Z_s\|^2\tr (dC_s)&=\|Z_t\|^2\tr (C_t)+\int_0^t\tr(C_s)\,d D_s\\
&-\int_0^t\tr(C_s)\,\tr(d K_s) -2\int_0^t\tr(C_s)\,d L_s\\
&=O(\|Z_t\|^2+D_t+\tr(K_t)+\tilde L_t).
\end{align*}
avec 
$$\tilde L_t:=\int_0^t\tr(C_s)\,d L_s.$$
Deux cas sont alors possibles \\\\
$\bullet\,$  soit $\langle \tilde L \rangle_{\infty}<\infty$ on déduit alors que  $\tilde L_t=O(D_t)$\\\\
$\bullet\,$  soit $\langle \tilde L \rangle_{\infty}=\infty$ et on conclut par la loi forte des 
 grands nombres pour les martingales scalaires que $\tilde L_t=o(\langle L\rangle_t).$
Ainsi on voit que 
$$\int_0^t \|Z_s\|^2\tr (dC_s)=O(\|Z_t\|^2+D_t+\tr(K_t)),$$
et par conséquent 
\begin{equation*}
\langle L\rangle_t=O\bl(\|Z_t\|^2+D_t +\tr(K_t)\br)\;\;p.s..
\end{equation*}
Encore une fois,  la loi forte des grands nombres pour les martingales scalaires  assure que 
\begin{equation}
 L_t=o\bl(\|Z_t\|^2+D_t+\tr(K_t)\br)\;\;p.s..
\end{equation}
Compte tenu du Lemme \ref{e5} et de la relation (\ref{e3}), on conclut que:
 \begin{equation}
L_t=o(A_t)\;\;p.s..
\end{equation}
\end{proof}
Par conséquent la relation 
\begin{equation}
A_t^{-1}\bigl(\|Z_t\|^2+\int_0^t {Z_s}^*SZ_s\,dA_s\bigr)\underset{t\rightarrow \infty}{\longrightarrow}
\tr\bigl(C^{1/2}SC^{1/2}\bigr) \;\;p.s.,
\end{equation}
découle aisément de la relation  
(\ref{e3}) et des lemmes \ref{e5}, \ref{e7} et \ref{e9}.

\vspace{.5cm}
Compte tenu de la preuve  du Théorème \ref{TLCPSQ} et sous les hypothèses  
$(\mathcal H1)$ et  $(\mathcal H)$, on a que 
$$\tilde{\mu}_t:=A_t^{-1}\int_0^t\delta_{Z_s}\,dA_s
\Longrightarrow \mu_\infty\;\;p.s..$$  
On en déduit alors que
$$\underset{t\rightarrow \infty}{{\underline{\lim}}} \int_{\mathbb R^d}x^*Sx \,d\tilde{\mu_t}(x) \geq 
\int_{{\mathbb R}^d} x^*Sx \,d\mu_{\infty}(x)\;\;p.s..$$ 
Or, d'après l'hypothèse { $(\mathcal H3)$}
$$\int_{\mathbb R^d} x^*Sx \, d\mu_{\infty}(x)=\tr\bigl(S\int x^*x \, d\mu_{\infty}(x)\bigr)=\tr(SC)
=\tr\bigl(C^{1/2}SC^{1/2}\bigr)$$ 
donc 
\begin{equation}\label{e13}
\underset{t\rightarrow \infty}{{\underline{\lim}}}A_t^{-1}\int_0^t {Z_s}^*SZ_s\,dA_s\geq\tr\bigl(C^{1/2}SC^{1/2}\bigr)\;\;p.s..
\end{equation}
Vu  les propriétés (\ref{11}) et (\ref{e13}), on conclut que 
\begin{equation}\label{e16}
\lim_{t\rightarrow \infty}A_t^{-1}\int_0^t {Z_s}^*SZ_s\,dA_s=\tr\bigl(C^{1/2}SC^{1/2}\bigr).
\end{equation}
 et on déduit  la loi du logarithme à savoir
$${ (L.L)} \;\;\;\|V_t^{-1}M_t\|^2=o(A_t)\;\;p.s.\;\;(t\rightarrow \infty).\hspace{6cm}$$
Par ailleurs, $S$ est inversible. 
On en déduit  alors, à l'aide de  
(\ref{e16}) et du fait que $$\log(\det V_t^2)\sim\tr(S)A_t\;\;(t\rightarrow\infty)
$$ (voir la relation (\ref{equiv})), la validité de la  propriété (LFQ). Ce qui achève la preuve.

\subsection{Preuve du Théorème \ref{TLCLFQ}.}
Posant 
$\theta_t :=V_t^{-1}\bigl(M_tM_t^*-[M]_t\bigr)(V_t^{*})^{-1},$
on a le lemme suivant
\begin{lemme}\label{ram}
\begin{multline*}
A_t^{-1/2} \int_0^t \tr(\theta_s)\,dA_s - 
\bigl(\log(\det V_t^2)\bigr)^{-1/2}\int_0^t\bigl(\tr[S]\bigr)^{-1/2}\tr(\theta_s)
d\bigl(\log(\det V_s^2)\bigr)\\\rightarrow 0\;\;{p.s.}\;\;\;\;(t\rightarrow \infty ).
\end{multline*}
\end{lemme}
\begin{proof}[Preuve.]
D'aprés la relation (\ref{equiv}) on voit que 
$$\Bigl[\frac{A_t}{\tr(S)\log(\det V_s^2)}\Bigr]^{-1/2} \rightarrow \tr(S)\;\;\;(t\rightarrow\infty).$$
D'autre part, en utlisant la relation (\ref{derive}) et la condition  ($\mathcal C3$) on obtient 
\begin{equation}
\lim_{t\rightarrow \infty}\frac{ \int_0^t\tr(\theta_s)\,dA_s}{\int_0^t\tr(\theta_s)
d\bigl(\log ( \det V_s^2 )\bigr)}=
\lim_{t\rightarrow\infty}\frac{a_t}{2\tr \bigl[V_t^{-1}\frac{dV_t}{dt}\bigr] }=\frac{1}{\tr(S) },
\end{equation}
ce qui achève la preuve du Lemme.
\end{proof}
Ainsi, on se ramène à démontrer que 
$$ A_t^{-1/2} \int_0^t \tr\bigl[\theta_s\bigr]\,dA_s \Rightarrow 2\sqrt{\tr\bigl(\tilde CRCR\bigr)}G$$
o\`u  $G$ est une gaussienne centrée réduite et $R$ désigne la matrice symétrique, positive solution de l'équation de Lyapounov:
$$I=RU+U^*R,$$
$U$ étant la matrice de la condition $(\mathcal C)$. 
\subsubsection{ Une relation fondamentale.}
Posant  $Z_t=V_t^{-1}M_t$, alors la formule d'It\^o donne:
\begin{multline}
\displaystyle Z_tZ_t^* =\int_0^tV_s^{-1}(dM_{s_-} )M_{s_-} ^*(V_s^{*})^{-1}+\int_0^tV_s^{-1}M_{s_-} (dM_{s_-} )^*(V_s^{*})^{-1}\\
                       -\int_0^tV_s^{-1}(dV_s)V_s^{-1}M_{s_-} M_{s_-} ^*(V_s^*)^{-1} 
                       +\int_0^tV_s^{-1}d[M]_s(V_s^*)^{-1}\hspace{0.9cm}\\
                       -\int_0^tV_s^{-1}M_{s_-} M_{s_-} ^*(V_s^{*})^{-1}(dV_s)^*(V_s^*)^{-1}.\hspace{4cm}
\end{multline}
Plus précisément, on a appliqué la formule d'Itô  à la  forme quadratique $(\langle u,Z_t\rangle^2)$ avec
$u\in \mathbb R^d$ et on obtient l'expression précédente par polarisation. 
En utilisant l'expression (\ref{IPP})
on obtient la relation fondamentale suivante:
\begin{equation}\label{rfond}
\theta_t + \int_0^tV_s^{-1}dV_s\theta_s+\int_0^t \theta_s(dV_s)^*(V_s^*)^{-1}
=H_t+H_t^*.
\end{equation}
avec
$$\displaystyle H_t=\int_0^tV_s^{-1}M_{s_-} (dM_{s_-} )^*(V_s^*)^{-1}.$$
Désormais, pour $u\in \mathbb R^d$, on pose:
\begin{equation}
\displaystyle H^u_t:=\int_0^tV_s^{-1}M_{s_-} (dM_{s_-} )^*(V_s^*)^{-1}u.
\end{equation}
Notre objectif est de démontrer un théorème limite centrale pour la martingale $H^u$. Pour celà, on va étudier
le comportement asymptotique du  crochet oblique de cette martingale ainsi que la condition 
de Lindeberg qui lui est associée. 


\subsubsection{Comportement asymptotique de  $(\langle H^u\rangle_t)_{t\geq 0}$}
Dans la suite, on démontre la proposition suivante
\begin{proposition}
\begin{equation}\label{crochet}
\frac{\langle H^u_t\rangle}{A_t}\rightarrow u^*\tilde Cu C\;\; \mbox{p.s..}
\end{equation}
\end{proposition}
\begin{proof}[Preuve.]
$\displaystyle (H^u_t)$ est une martingale vectorielle de variation quadratique prévisible:
\begin{multline}\label{Hu}
\langle H^u\rangle_t =\int_0^tV^{-1}_sM_{s_-} \bigl[u^*V_s^{-1}d\langle M\rangle_s(V_s^*)^{-1}u\bigr]
(M_{s_-} )^*(V_s^{*})^{-1}\\
                     =\int_0^t Z_s\bigl[u^*V_s^{-1}d\langle M \rangle_s(V_s^*)^{-1}u\bigr]Z_s^*.\hspace{4.1cm}
\end{multline}
On en  déduit alors que pour tout $x\in \mathbb R^d$ on a:
\begin{align*}
x^*\langle H^u\rangle_tx&=\int_0^t x^*Z_s\bigl[u^*V_s^{-1}d\langle M\rangle_s(V_s^*)^{-1}u\bigr]Z_s^*x\\
                        &=\int_0^t\langle x,Z_s\rangle^2\bigl[u^*V_s^{-1}d\langle M \rangle_s(V_s^*)^{-1}u\bigr].
\end{align*}
On pose alors $$F_t(u):=\int_0^t\exp(A_s)u^*V^{-1}_sd\langle M\rangle_s(V^*_s)^{-1}uds.$$ 
Ainsi par la formule d'intégration par partie on déduit que:
\begin{align}\label{Gu}
x^*\langle H^u\rangle_t x&=\exp(-A_t)F_t(u)\langle x,Z_t\rangle^2+\int_0^t
\exp(-A_s)F_s(u)\langle x,Z_s\rangle^2d(A_s)-G_t
\end{align}
avec $$G_t:=\int_0^t\exp(-A_s)F_s(u)d(\langle x,Z_s\rangle^2).$$
Le résultat suivant est utile pour la suite
\begin{lemme}\label{F}
\begin{equation}
\exp(-A_t)F_t(u)\rightarrow u^*\tilde Cu\;\;\;\mbox{p.s.,}
\end{equation}
avec $\tilde C=UC+CU^*$; $U$   étant la matrice introduite dans $\mathcal (C3)$. 
\end{lemme}
\begin{proof}[Preuve.]
Par la formule d'intégration par parties, on obtient:
\begin{align*}
F_t(u)&=\exp(A_t)u^*V_t^{-1}\langle M\rangle_t(V_t^*)^{-1}u-\int_0^tu^*V_s^{-1}
\langle M\rangle_s(V_s^*)^{-1}ud\bigl(\exp(A_s)\bigr)\\
      &+\int_0^tu^*\bigl(a_s^{-1}V_s^{-1}\frac{dV_s}{ds}\bigr)V_s^{-1}
\langle M\rangle_s(V_s^*)^{-1}ud\bigl(\exp(A_s)\bigr)\\
      &+\int_0^tu^*V^{-1}_s\langle M\rangle_s(V_s^*)^{-1}\bigl(a_s^{-1}V_s^{-1}
\frac{dV_s}{ds}\bigr)^*ud\bigl(\exp(A_s)\bigr).
\end{align*}
D'après l'hypothèse  $(\mathcal H1)$ et le lemme de Toeplitz on voit que:
\begin{multline*}
\frac{1}{\exp(A_t)-1}\Bigl(\exp(A_t)u^*V_t^{-1}\langle M\rangle_t(V_t^*)^{-1}u\\
-\int_0^tu^*V_s^{-1}\langle M\rangle_s(V_s^*)^{-1}ud\bigl(\exp(A_s)\Bigr)\rightarrow 0\;\;\;\mbox{p.s.}
\end{multline*}
et d' autre part  d'après l'hypothèse $(\mathcal H1)$, 
le lemme de Toeplitz et la condition  $\mathcal (C3)$ on obtient que:
\begin{multline*}
\frac{1}{exp(A_t)-1}\Bigl(\int_0^tu^*\bigl(a_s^{-1}V_s^{-1}\frac{dV_s}{ds}\bigr)V_s^{-1}
\langle M\rangle_s(V_s^*)^{-1}ud\bigl(\exp(A_s)\bigr)\\
      +\int_0^tu^*V^{-1}_s\langle M\rangle_s(V_s^*)^{-1}
\bigl(a_s^{-1}V_s^{-1}\frac{dV_s}{ds}\bigr)^*ud\bigl(\exp(A_s)\bigr)\Bigr)
\rightarrow \tilde C=UC+CU^*\;\;\;\mbox{p.s.}
\end{multline*}
ce qui achève la démonstration.
\end{proof}
Par conséquent, en combinant le lemme précédent et
 la propriété { (LL)} on voit immédiatement que:
\begin{equation}\label{H1} 
\exp(-A_t)F_t(u)\langle x,Z_t\rangle^2=o(A_t)\;\;\;\mbox{p.s.}
\end{equation}
Par ailleurs, d'aprés  le lemme \ref{F} et la propriété { (LFQ)}, il vient que:
\begin{equation}\label{H2}
A_t^{-1}\int_0^t\exp(-A_s)F_s(u)\langle x,Z_s\rangle^2d(A_s)\rightarrow u^*\tilde Cu x^*Cx\;\;\;\mbox{p.s..}
\end{equation}
Dans la suite on s'interresse au comportement asymptotique de $G$.
Vu le Lemme \ref{F}, il est clair que 
\begin{equation}
G_t=O\bigl(\langle x,Z_t\rangle^2\bigr)
\end{equation}
donc par la propriété (LL) on a
\begin{equation}\label{H3}
G_t=o(A_t)\;\;\; {p.s..}
\end{equation}
Compte tenu de l'expression (\ref{Gu}) des relations (\ref{H1}),(\ref{H2}) et (\ref{H3})
on conclut que 
\begin{equation}
\frac{\langle H^u_t\rangle}{A_t}\rightarrow u^*\tilde Cu C\;\; \mbox{p.s..}
\end{equation}
Ce qui achève la preuve.
\end{proof}
 
\subsubsection{Vérification de la condition de Lindeberg pour la martingale $( H^u_t)_{t\geq 0}$}		

\begin{definition}
Soient $A=(A_{t}),B=(B_{t})$ deux processus croissants issus de 0. On dit
que $A$ est \textbf{domin\'{e} au sens fort} par $B,$et on \'{e}crit: $A<<B,$
si $(B_{t}-A_{t};t\geq 0)$ est un procssus croissant.
\end{definition}
Le r\'{e}sultat utile suivant est \'{e}vident:
\begin{lemme}
Si $A<<B$, leurs compensateurs pr\'{e}visibles $\widetilde{A%
},$ $\widetilde{B}$ v\'{e}rifient aussi $\widetilde{A}<<\widetilde{B\text{ }}.$
\end{lemme}
\textbf{Application \`{a} la martingale }$(H_{t})$\\
Le saut \`{a} l'instant $t$ de la martingale matricielle:
\begin{equation*}
H_{t}=\int_{0}^{t}Z_{s-}d(^{\ast }M_{s})^{\ast }V_{s}^{-1},\text{ }%
Z_{t}=V_{t}^{-1}M_{t}\text{ },
\end{equation*}
vaut: 
\begin{equation*}
\Delta H_{t}=Z_{t-}{}^{\ast }(\Delta M_{t})^{\ast }V_{t}^{-1};
\end{equation*}
donc: 
\begin{align*}
\|\Delta H_{t}\|^{2} &=tr\{\Delta H_{t}\text{ }%
^{\ast }\Delta H_{t}\}=\|Z_{t-}\|^{2}
\| V_{t}^{-1}\Delta M_{t}\text{ }\|^{2} \\
&=\|Z_{t-}\|^{2}V_{t}^{-1}\Delta \lbrack
M]_{t}^{\ast }V_{t}^{-1}=\| Z_{t-}\|^{2}\Delta
\Lambda _{t}
\end{align*}
o\`{u} $(\Lambda _{t})$ est le processus croissant: 
\begin{equation*}
\Lambda _{t}=\int_{0}^{t}V_{s}^{-1}d[M]_{s}^{\ast }V_{s}^{-1}.
\end{equation*}
Pour $r>0,$ $t>0,$ posant: 
\begin{equation*}
\sigma _{t}^{H}(r)=\sum_{s\text{ }\leq \text{ }t}\| \Delta
H_{s}\| ^{2}\mathbf{1}_{\{\text{ }\| \Delta
H_{s}\| >r\text{ }\}},
\end{equation*}
la condition de Lindeberg au sens de la convergence presque sure pour la
martingale $H$ s'\'{e}crit: 
\begin{equation}\label{lind}
\forall \epsilon >0,\text{ \ }A_{t}^{-1}\widetilde{\sigma _{t}^{H}(\epsilon 
\sqrt{A_{t}})}\rightarrow 0,\text{ \ }t\rightarrow +\infty ,\text{ p.s. }
\end{equation}
Pour \'{e}tablir ce r\'{e}sultat, on exploite les deux lemmes suivants:
\bigskip
\begin{lemme}
L'hypoth\`{e}se ($\mathcal H2$) implique que presque surement 
$$\sup_{t\text{ }>\text{ }0}\| V_{t}^{-1}\Delta
M_{t}\| <+\infty. $$
\end{lemme}
\begin{proof}[Preuve]
En effet, on a: 
\begin{equation*}
\sum_{s\text{ }\leq \text{ }. }\Delta M_{s}\text{ }^{\ast }(\Delta
M_{s})<<[M]_{. }
\end{equation*}

ce qui implique que 
\begin{equation*}
\sum_{s\text{ }\leq \text{ }t}\| V_{t}^{-1}\Delta M_{s}\text{ }%
\| ^{2}\leq tr\{V_{t}^{-1}[M]_{t}^{\ast }V_{t}^{-1}\};
\end{equation*}

d'o\`{u} le r\'{e}sultat du lemme, car 
\begin{equation*}
\| V_{t}^{-1}\Delta M_{t}\text{ }\| ^{2}\leq \sum_{s%
\text{ }\leq \text{ }t}\| V_{t}^{-1}\Delta M_{s}\text{ }%
\| ^{2}\leq tr\{V_{t}^{-1}[M]_{t}^{\ast }V_{t}^{-1}\}=O(1)\text{
p.s. .}
\end{equation*}
\end{proof}
\begin{lemme}
Etant donn\'{e} $\alpha \in ]0,1],$ alors: 
\begin{equation*}
\sigma _{. }^{H}(\alpha ^{-3})<<\sigma _{. }^{1}(\alpha
^{-1})+\sigma _{. }^{2}(\alpha ^{-1})
\end{equation*}
 o\`{u} pour $t>0,$ $r>0:$%
\bigskip
\begin{equation*}
\sigma _{t}^{1}(r)=\sum_{s\text{ }\leq \text{ }t}\| \Delta
H_{s}\| ^{2}\mathbf{1}_{\{\text{ }\|
Z_{s-}\| >r\text{ }\}}\text{ \ \ };\text{ \ }\sigma
_{t}^{2}(r)=\sum_{s\text{ }\leq \text{ }t}\| \Delta
H_{s}\| ^{2}\mathbf{1}_{\{\text{ }\|
V_{s}^{-1}\Delta M_{s}\| >r\text{ }\}}
\end{equation*}
\end{lemme}
\begin{proof}[Preuve.]
L'assertion du lemme d\'{e}coule de la d\'{e}composition 
\'{e}vidente suivante: 
\begin{align*}
\sigma _{t}^{H}(\alpha ^{-3}) &=\sum_{s\text{ }\leq \text{ }%
t}\| \Delta H_{s}\| ^{2}\mathbf{1}_{\{\text{ }%
\| \Delta H_{s}\| >\alpha ^{-3}\text{ }%
,\| Z_{s-}\| >\alpha ^{-1}\text{ }\}}\\&+\sum_{s\text{ }%
\leq \text{ }t}\| \Delta H_{s}\| ^{2}\mathbf{1}_{\{%
\text{ }\| \Delta H_{s}\| >\alpha ^{-3}\text{ },%
\text{ }\| Z_{s-}\| \leq \text{ }\alpha ^{-1},\text{ 
}\| V_{s}^{-1}\Delta M_{s}\| \leq \text{ \ }\alpha
^{-1}\text{ }\}} 
\\&+\sum_{s\text{ }\leq \text{ }t} \| \Delta H_{s}\|
^{2}\mathbf{1}_{\{\text{ }\| \Delta H_{s}\| >\alpha
^{-3}\text{ },\text{ }\| Z_{s-}\| \leq \text{ }%
\alpha ^{-1},\text{ }\| V_{s}^{-1}\Delta M_{s}\| >%
\text{\ }\alpha ^{-1}\text{ }\}}.
\end{align*}
Il est clair que le premier (resp le troisi\`{e}me) terme du membre de
droite de cette \'{e}galit\'{e} est domin\'{e} au sens fort par $(\sigma
_{t}^{1}(\alpha ^{-1}))$ \ (resp. $(\sigma _{t}^{2}(\alpha ^{-1}))=(\sigma
_{t}^{2}(\min (\alpha ^{-1},\alpha ^{-2}))$ ). Pour le deuxi\`{e}me terme,
on remarque 
\begin{multline*}
\{\text{ }\| \Delta H_{s}\| >\alpha ^{-3},\text{ }%
\| Z_{s-}\| \leq \alpha ^{-1},\text{ }\|
V_{s}^{-1}\Delta M_{s}\| \leq \ \alpha ^{-1}\}\subset \\
\{\text{ }%
\| Z_{s-}\| \text{ }\leq \alpha ^{-1},\text{ }%
\| Z_{s-}\| \geq \alpha ^{-2}\text{ }\}=\emptyset 
\text{ p.s.}
\end{multline*}
d'apr\`{e}s le choix de $\alpha .$ Le lemme est \'{e}tabli.
\end{proof}

\begin{corollaire}
Avec les notations du lemme 3, on a: 
\begin{equation*}
\sigma _{. }^{1}(r)<<\int_{0}^{. }\|
Z_{s-}\| ^{2}\mathbf{1}_{\{\text{ }\|
Z_{s-}\| >r\text{ }\}}d\Lambda _{s}\text{ \ , \ \ }\sigma
_{. }^{2}(r)<<\int_{0}^{. }\| Z_{s-}\|
^{2}\mathbf{1}_{\{\Delta \text{ }\Lambda _{s}\text{ }>\text{ }r^{2}\text{ }%
\}}d\Lambda _{s}.
\end{equation*}

Par cons\'{e}quent, pour tout $t\geq 0$: 
\begin{eqnarray*}
\widetilde{\sigma _{t}^{1}(r)}\text{ } &\leq &\int_{0}^{t}\|
Z_{s-}\| ^{2}\mathbf{1}_{\{\text{ }\|
Z_{s-}\| >r\text{ }\}}d\widetilde{\Lambda _{s}}\text{ , \ } \\
\widetilde{\sigma _{t}^{2}(r)}\text{ } &\leq &\int_{0}^{t}\|
Z_{s-}\| ^{2}\mathbf{1}_{\{\text{ }\Delta \Lambda _{s}>r^{2}\text{
}\}}d\widetilde{\Lambda _{s}}\leq \mathbf{1}_{\{\text{ }\sup_{s\text{ }\leq 
\text{ }t}\Delta \Lambda _{s}\text{ }>\text{ }r^{2}\text{ }%
\}}\int_{0}^{t}\| Z_{s-}\| ^{2}d\widetilde{\Lambda
_{s}}\text{ },
\end{eqnarray*}

avec $$\widetilde{\Lambda _{t}}=\int_{0}^{t}V_{s}^{-1}d<M>_{s}^{\ast
}V_{s}^{-1}.$$
\end{corollaire}
\textbf{ Validit\'{e} de la condition de Lindeberg}\\

Compte tenu de ce qui pr\'{e}c\`{e}de, pour tous $\alpha \in ]0,1[$ , $%
\epsilon >0$ et $t>0,$ on a 
\begin{equation*}
\overline{\lim_{t}}\text{\ }A_{t}^{-1}\widetilde{\sigma _{t}^{H}(\epsilon 
\sqrt{A_{t}})}\text{ }\leq \overline{\lim_{t}}\text{\ }A_{t}^{-1}\widetilde{%
\sigma _{t}^{H}(\alpha ^{-3})}.
\end{equation*}
Mais
\begin{multline*}
\widetilde{\sigma _{t}^{H}(\alpha ^{-3})}\text{ }\leq
\int_{0}^{t}\| Z_{s-}\| ^{2}\mathbf{1}_{\{\text{ }%
\| Z_{s-}\| >\alpha ^{-1}\text{ }\}}d\widetilde{%
\Lambda _{s}}\text{ }\\+\mathbf{1}_{\{\text{ }\sup_{s\text{ }\leq \text{ }%
t}\Delta \Lambda _{s}\text{ }>\text{ }\alpha ^{-2}\text{ }\}}\int_{0}^{t}%
\| Z_{s-}\| ^{2}d\widetilde{\Lambda _{s}},
\end{multline*}

donc par la loi forte quadratique, pour tous $\alpha \in ]0,1[$ , $\epsilon
>0$: 
\begin{multline*}
\overline{\lim_{t}}\text{\ }A_{t}^{-1}\widetilde{\sigma _{t}^{H}(\epsilon 
\sqrt{A_{t}})}\text{ }\leq \int_{0}^{+\infty }\| x\|
^{2}\mathbf{1}_{\{\text{ }\| x\| >\alpha ^{-1}\text{ 
}\}}d\mu _{\infty }(x)
\\+\mathbf{1}_{\{\text{ }\sup_{s\text{ }\geq \text{ }%
0}\Delta \Lambda _{s}\text{ }>\text{ }\alpha ^{-2}\text{ }\}}\text{ }%
\int_{0}^{+\infty }\| x\| ^{2}d\mu _{\infty }(x)%
\text{ \ p.s.}
\end{multline*}

ce qui implique que la condition de Lindeberg est v\'{e}rifi\'{e}e.

\bigskip



\subsubsection{Fin de la preuve du Théorème \ref{TLCLFQ}}
Vu les propriétés  (\ref{crochet}) et  (\ref{lind}),  le (TLCG) s'applique pour la martingale vectorielle $H^u$
et on a: 
\begin{equation*}
A_t^{-1/2}H_t\Rightarrow \mathfrak N_{\s{d\times d}}(0,\,\tilde C\otimes C)
\end{equation*}
o\`u $\tilde C=(UC+CU^*)$ . De l'égalité (\ref{rfond}), il vient  que pour toute matrice symétrique positive $R$ on a
\begin{multline*}
 A_t^{-1/2}\tr\bigl(R\theta_t \bigr)+ A_t^{-1/2}\int_0^t \tr\Bigl\{\bigl[RV_s^{-1}dV_s +(V_s^{-1}dV_s)^*R\bigr]  \theta_{s_-} \Bigr\}
\\\Rightarrow 2\sqrt{\tr\bigl(\tilde CRCR\bigr)}G
\end{multline*}
o\`u  $\theta_t =V_t^{-1}\bigl(M_tM_t^*-[M]_t\bigr)(V_t^{*})^{-1}$ et $G$  une variable aléatoire gaussienne
 centrée réduite. Par conséquent, comme $Z_t=V_t^{-1}M_t$ converge en loi et comme  $V_t^{-1}[M]_t(V_t^{*})^{-1}$
 converge p.s., on en déduit que 
\begin{multline*}
 A_t^{-1/2}\int_0^t \tr\Bigl\{\Bigl[R\bigl(a_s^{-1}V_s^{-1}\frac{dV_s}{ds}\bigr) 
+\bigl(a_s^{-1}V_s^{-1}\frac{dV_s}{ds}\bigr)^*R\Bigr]  \theta_s \Bigr\}\,dA_s\\ 
\Rightarrow 2\sqrt{\tr\bigl(\tilde CRCR\bigr)}G
\end{multline*}
Compte tenu de la condition { $(\mathcal C3$)} et de la ${(L.L)}$ qui garantie que $\theta_t =o(A_t)\;\;p.s.$,
 on obtient que 
\begin{multline}\label{negliger}
\int_0^t \tr\Bigl\{\Bigl[R\bigl(a_s^{-1}V_s^{-1}\frac{dV_s}{ds}-U\bigr) +
\bigl(a_s^{-1}V_s^{-1}\frac{dV_s}{ds}-U\bigr)^*R\Bigr]  \theta_s \Bigr\}\,dA_s\\=
o\Bigl(\int_0^tA_s^{-1/2}\,dA_s\Bigr)=o(A_t^{1/2})\;\;p.s..
\end{multline}
Ainsi pour $R$ solution de l'équation de Lyapounov $I=RU+U^*R$ il vient que:
\begin{equation}
 A_t^{-1/2}\int_0^t \tr\bigl\{\theta_s\bigr\}\,dA_s\\\Rightarrow 
2\sqrt{\tr\bigl(\tilde CRCR\bigr)}G
\end{equation}
On conclut la preuve de la première partie du théorème par le Lemme \ref{ram}. La deuxième partie du théorème est 
immédiate en remarquant que l'hypothèse ajoutée est équivalente à:\\

\hfill 
$A_t^{^{\rho}}\Bigl|\tr\bigl\{V_t^{-1}\bigl([M]_t\bigr)(V_t^{*})^{-1}-C\bigr\}\Bigr|=O(1)\;\;\mbox{ p.s.},\;\;\rho>1/2.$
\hfill $\Box$
\subsection{Preuve du Théorème \ref{LLIL}}
D'après la relation (\ref{rfond}) et pour toute matrice symétrique, solution de l'équation de Lyapounov $RU+U^*R=I$, on 
a:
\begin{equation}
\tr\bigl(R\theta_t \bigr)+ \int_0^t \tr\Bigl\{\Bigl[R\bigl(a_s^{-1}V_s^{-1}\frac{dV_s}{ds}\bigr) 
+\bigl(a_s^{-1}V_s^{-1}\frac{dV_s}{ds}\bigr)^*R\Bigr]  \theta_s \Bigr\}\,dA_s =2\tr\{RH_t
\}.
\end{equation}\\
Comme,  $$\;\mathbb E\bl[\sup_t M_t(\Delta M_t)^*\br]< \infty\;\;\;\mbox{implique que}\; \;\;\mathbb E\bl[\sup_t\Delta \tr\{RH_t\}
\br]
<\infty,$$ on en déduit par le Théorème 3 de \cite{lepingle} et par la relation (\ref{crochet}) que la martingale scalaire $(\tr\{RH_t\})_{t\geq 0}$ vérifie
 une loi du logarithme itéré donnée par:
\begin{equation}
{{\underset{t\rightarrow \infty}{\overline{\lim}}}}\;\frac{\tr\{RH_t\}}{h(A_t)}\leq \sqrt{\tr\bigl(\tilde CRCR\bigr)}\;\;\;\mbox{p.s.,}
\end{equation}
o\`u $h(u)=\sqrt{2u\log\log\, u}\;\;\mbox{pour\;} u\geq e$.
Compte tenu de la relation (\ref{negliger}) et du fait que $\theta_t=o(A_t)$ on obtient que:
\begin{equation*}
{{\underset{t\rightarrow \infty}{\overline{\lim}}}}\;\frac{1}{h(A_t)}\int_0^t \tr\Bigl\{ V_s^{-1}\bigl(M_sM_s^*-[M]_s\bigr)(V_s^{*})^{-1}  \Bigr\}\,dA_s \leq \sqrt{\tr\bigl(\tilde CRCR\bigr)}\;\;\;\mbox{p.s..}
\end{equation*}
La fin de la preuve est similaire à celle de la preuve précédente..
\hfill $\Box$
\subsection{Preuve des Théorèmes  \ref{tlclfqs} et \ref{llils}}
Comme La normalisation scalaire est un cas particulier de la normalisation matricielle alors
de l'égalité (\ref{rfond}) et de la condition (\ref{speed}) on voit que  
$$ A_t^{-1/2}\theta_t + 2A_t^{-1/2}\int_0^t  \theta_{s_-}v_s^{-1}dv_s \Rightarrow 2\sqrt{2\eta\, }C G$$
o\`u  $\theta_t =v_t^{-2}\bigl(M_tM_t^*-[M]_t\bigr)$ et $G$  une variable aléatoire gaussienne
 centrée réduite. Par conséquent, comme $Z_t=v_t^{-1}M_t$ converge en loi et comme  $v_t^{-2}[M]_t$
 converge p.s., on en déduit que 

$$ A_t^{-1/2}\int_0^t  \theta_{s_-} a_s^{-1}v_s^{-1}v'_s\,dA_s \Rightarrow \sqrt{2\eta\, }C G$$

Compte tenu de la condition { (\ref{speed})} et de la ${(L.L)}$ qui garantie que $\theta_t =o(A_t)\;\;p.s.$,
 on obtient 
\begin{equation}\label{negligers}
\int_0^t \bigl( a_s^{-1}v_s^{-1}v'_s-\eta \bigr)\,dA_s=
o\Bigl(\int_0^tA_s^{-1/2}\,dA_s\Bigr)=o(A_t^{1/2})\;\;p.s..
\end{equation}
Ainsi on en déduit que 
\begin{equation}
 A_t^{-1/2}\int_0^t \tr\bigl\{\theta_s\bigr\}\,dA_s\\\Rightarrow 
(\sqrt{2\eta}\,C)\,G
\end{equation}
On conclut la preuve de la première partie du théorème par le Lemme \ref{ram} qui reste valable pour une
normalisation $V_t$ scalaire. La deuxième partie du Théorème \ref{tlclfqs} est 
immédiate en remarquant que l'hypothèse ajoutée est équivalente à:
$$A_t^{^{\rho}}\Bigl|\bigl\{v_t^{-2}\bigl([M]_t\bigr)-C\bigr\}\Bigr|=O(1)\;\;\mbox{ p.s.},\;\;\rho>1/2.$$\\

On donne à présent la preuve du  Théorème \ref{llils}. En effet,   
de la relation (\ref{rfond}) on 
a:
\begin{equation}
 \theta_t + 2\int_0^t \Bigl\{\bigl(a_s^{-1}v_s^{-1}{v'_s}\bigr) \theta_s \Bigr\}\,dA_s =2 H_t.
\end{equation}\\
Comme,  $$\;\mathbb E\bl[\sup_t M_t(\Delta M_t)^*\br]< \infty\;\;\;\mbox{implique que}\; \;\;
\mathbb E\bl[\sup_t\Delta H_t
\br]
<\infty,$$ on en déduit par le Théorème 3 de \cite{lepingle} et par la relation (\ref{crochet}) 
(valable pour une normalisation scalaire) que la martingale  $(H_t)_{t\geq 0}$ vérifie
 une loi du logarithme itéré donnée par:
\begin{equation}
{{\underset{t\rightarrow \infty}{\overline{\lim}}}}\;\frac{H_t}{h(A_t)}\leq 
\sqrt{ 2\eta}\,C\;\;\;\mbox{p.s.,}
\end{equation}
o\`u $h(u)=\sqrt{2u\log\log\, u}\;\;\mbox{pour\;} u\geq e$.
Compte tenu de la relation (\ref{negliger}) et du fait que $\theta_t=o(A_t)\;\;p.s.$ on obtient que:
\begin{equation*}
{{\underset{t\rightarrow \infty}{\overline{\lim}}}}\;\frac{1}{h(A_t)}\int_0^t 
v_s^{-2}\bigl(M_sM_s^*-[M]_s\bigr) \,dA_s \leq \sqrt{2\eta }\,C\;\;\;\mbox{p.s..}
\end{equation*}
La fin de la preuve est similaire à celle de la preuve précédente..
\hfill $\Box$
\section{Application: Estimation de la variance d'un \\P.A.I.S. pondéré}
Une question interessante nous  a été posé au fur et à mesure que 
ce travail progressait. En effet, il s'agissait de savoir si on pouvait améliorer la
vitesse logarithmique (lente)  dans la propriété (TLCPS) ainsi que dans les autres propriétés
qui lui sont associées. Dans ce qui suit on donne une réponse à cette question sous forme d'application.
 On se propose alors d'estimer
la variance d'un P.A.I.S. pondéré. On dira que le processus
 $(\tilde S_t)_{t\geq 0}$ est un P.A.I.S. pondéré  s'il est de la forme
$$\tilde S_t:= \int_0^tw_s\,dS_s$$
o\`u $w$ est un processus à variation fini alors que
$S$ est un processus à accroissement indépendants et stationnaires dont la mesure de
Lévy des sauts $\nu$ vérifie:
\begin{equation}\label{pais}
\nu(dt,dx)=dt\,F(dx),\;\;\mbox{avec}\;\;\int|x|^{2p}F(dx)<\infty\;\;\mbox{pour un}\;\;  p>1,
\end{equation}
o\`u $F$ est une mesure positive sur $\mathbb R$. 
On note:
$$m=\mathbb E\,S_1 ,\;\;\;\;\sigma^2=\mathbb E\, S_1^2-m^2\;\;\;\mbox{et}\;\;\;\tilde N_t=\int_0^tw_r\,d(S_r-mr).$$
\begin{proposition}\label{appli}
Avec les notations précédentes et pour $$w_t=\frac{t^{-\alpha /2}}{1-\alpha}\exp{\frac{t^{1-\alpha}}{2(1-\alpha)}},\;\;
\;\;\alpha\in(0,1)$$
on a les propriétés\\
\begin{enumerate}
\item (TLCPS)
\begin{center}
$\;\;\displaystyle \frac{1-\alpha}{t^{1-\alpha}} \int_0^t
\delta_{\Bigl\{e^{-\frac{s^{1-\alpha}}{2(1-\alpha)}}\tilde N_s\Bigr\}}\,\frac{ds}{s^{\alpha}}
\Rightarrow \mathfrak N(0,\sigma^2)$
\end{center}
\item (LFQ)
\begin{center}
$\;\;\displaystyle \tilde{\sigma_t}:=\frac{1-\alpha}{t^{1-\alpha}} \int_0^t  e^{-\frac{s^{1-\alpha}}{1-\alpha}}\,
\tilde N_s^2\,\frac{ds}{s^{\alpha}}
\rightarrow \sigma^2\;\;p.s.\;\;(t\rightarrow \infty).$
\end{center}
\vspace{.5cm}
Si de plus pour $\rho>1/2$  on a
$$ \exp{\Bigl\{-\frac{t^{1-\alpha}}{2(1-\alpha)}\Bigr\}}\sum_{s\leq t}(\Delta\tilde S_s)^2-\int_Rx^2F(dx)=
O(t^{\rho(1-\alpha)})\;\;(t\rightarrow \infty)
$$
alors on a la propriété 
\item (TLCLFQ)
\begin{center}
$\;\;t^{\frac{1-\alpha}{2}}(\tilde{\sigma_t}-\sigma^2)\Rightarrow \mathfrak N\bigl(0,4(1-\alpha)\sigma^4\bigr).$
\end{center}
\end{enumerate}
\end{proposition}
La preuve de la proposition est laissée en annexe.\\
{\bf Remarques}
\begin{enumerate}
\item Les preuves des propriétés donées dans le paragraphe 1.2 et celle de la proposition précédente sont similaires. 
\item Dans la proposition \ref{appli}, de la propriété (LFQ)   on voit que   $\tilde \sigma_t$ 
est un estimateur fortement consistant
de $\sigma^2$. Cependant, vu la proprété (TLCLFQ), l'intervalle de confiance associé à cet estimateur est asymptotiquement
meilleur que celui donné par l'estimateur sans pondération à savoir l'estimateur
 $\hat\sigma$ (voir partie 1.2)
\end{enumerate} 

\section{Annexe}
\subsection{Preuve du Lemme \ref{chaab3}}
La propri\'{e}t\'{e}  (\ref{EQ1}) implique que pour tout couple $\left( \rho
_{1},\rho _{2}\right) \in \Bbb{R}_{+}\times \Bbb{R}_{+}$ avec 
$\rho _{1}\leq \rho _{2},$ on a:
\begin{equation}\label{EQ2}
\tr\Bigl(
\int\nolimits_{\rho_1}^{\rho_2}V_{s}^{-1}d\left( V_{s}V^*_{s}\right)
(V_{s}^*)^{-1}\Bigr)  =\log (\det V_{\rho _{2}})^{2}-\log (\det V_{\rho
_{1}} )^{2}.
\end{equation}
Or, vu que:
\begin{align*}
\tr\Bigl(
\int\nolimits_{\rho_1}^{\rho_2}V_{s}^{-1}d\left( V_{s}V^*_{s}\right)
(V_{s}^*)^{-1}\Bigr)&= \tr
\Bigl( \int\nolimits_{\rho _{1}}^{\rho _{2}}\left( V_{s}V_{s}^*\right) ^{-1}d\left( V_{s}V_{s}^*\right) \Bigr)\\
&\geq \tr\Bigl( \left( V_{\rho _{2}}V^*_{\rho _{2}}\right) ^{-1}\int\nolimits_{\rho _{1}}^{\rho _{2}}d\left(
V_{s}V^*_{s}\right) \Bigr)\\
&\geq \,\tr\Bigl(\left( V_{\rho _{2}}V^*_{\rho _{2}}\right) ^{-1}\left( V_{\rho _{2}}V^*_{\rho_{2}}-
V_{\rho _{1}}V^*_{\rho _{1}}\right) \Bigr)\\ 
&\geq \,\tr\Bigl( I_d-V_{\rho _{2}}^{-1}V_{\rho_{1}}V^*_{\rho _{1}}(V^*_{\rho _{2}})^{-1}\Bigr), 
\end{align*}
on en d\'{e}duit que:
\begin{equation}\label{EQ3}
\tr\Bigl(
\int\nolimits_{\rho_1}^{\rho_2}V_{s}^{-1}d\left( V_{s}V^*_{s}\right)
(V_{s}^*)^{-1}\Bigr)\geq d-\left\| V_{\rho _{2}}^{-1}V_{\rho _{1}}\right\|
^{2}.
\end{equation}
Les deux propri\'{e}t\'{e}s (\ref{EQ2}) et (\ref{EQ3}) impliquent donc que:
\begin{equation}\label{EQ4}
d-\left\| V_{\rho _{2}}^{-1}V_{\rho
_{1}}\right\| ^{2}\leq  \log (\det V_{\rho _{2}})^{2}-\log (\det V_{\rho
_{1}})^{2}.
\end{equation}
Pour un $n_0$ fixé consid\'{e}rons  maitenant la subdivision suivante: 
$\rho _{0}=\rho <\rho
_{1}<\cdots <\rho _{n_0}=r$, on a alors:
\begin{align*}
\left\| V_{r}^{-1}V_{\rho }\right\| ^{2}&\leq
\prod\limits_{j=0}^{n_0-1}\left\| V_{\rho _{j+1}}^{-1}V_{\rho _{j}}\right\| 
^{2}\\&=\prod\limits_{j=0}^{n_0-1}\left[ d-\left( d-\left\| V_{\rho
_{j+1}}^{-1}V_{\rho _{j}}\right\| ^{2}\right) \right] \\
&=d^{n_0}\prod\limits_{j=0}^{n_0-1}\left[ 1-\left( 1-\frac{\left\| V_{\rho
_{j+1}}^{-1}V_{\rho _{j}}\right\| ^{2}}{d}\right) \right] \\
&\leq d^{n_0}\exp \left\{ -\sum\limits_{j=0}^{n-1}\left( 1-
\frac{\left\| V_{\rho
_{j+1}}^{-1}V_{\rho _{j}}\right\| ^{2}}{d}\right) \right\}\\ 
&\leq d^{n_0}\exp \left\{ -\frac{1}{d}\sum\limits_{j=0}^{n-1}\left[ \log (\det V_{\rho
_{j+1}})^{2}-\log (\det V_{\rho _{j}})^{2}\right] \right\} .
\end{align*}
D'o\`{u} l'in\'{e}galit\'{e}:\\

\hfill
$\displaystyle\left\| V_{r}^{-1}V_{\rho }\right\| ^{2}\leq d^{n_0} 
\Bigl(\dfrac{\det V_{\rho }}{\det V_{r}}\Bigr)^{\frac{2}{d}}.$\hfill $\Box$
\subsection{Preuve de la Proposition \ref{appli}}
On sait que $S$ est un processus à accroissement indépendents et stationnaires (P.A.I.S.) par conséquent 
$\tilde N_t:=\int_0^tw_r\,d(S_r-mr)$ est une martingale dont la variation quadratique  est donnée par
$\langle N \rangle_t=\sigma^2\int_0^tw_r^2\,dr$. Pour 
$$w_t=\frac{t^{-\alpha /2}}{1-\alpha}\exp{\frac{t^{1-\alpha}}{2(1-\alpha)}},\;\;
\;\;\alpha\in(0,1)$$
on voit que $$e^{\frac{t^{1-\alpha}}{1-\alpha}}\langle N \rangle_t\rightarrow \sigma^2,\;\;\;(t\rightarrow \infty).$$
Ainsi l'hypothèse ($\mathcal H1$) est vérifiée.
L'hypothèse ($\mathcal H'2$) est immédiate. En effet
\begin{align*}
\int_{\mathbb R}\int_0^te^{-\frac{t^{1-\alpha}}{1-\alpha}}|x|^2\nu^{\tilde M}(ds,dx)=&
e^{-\frac{t^{1-\alpha}}{1-\alpha}}\int_0^t\frac{s^{-\alpha /2}}{1-\alpha}e^{\frac{s^{1-\alpha}}{2(1-\alpha)}}\,ds
\int_{\mathbb R }|x|^2F(dx)\\&\,\rightarrow 0\;\;\;(t\rightarrow \infty).
\end{align*}
La propriété (TLCPS) est prouvée. Afin de démontrer les propriétés (LFQ) et (TLCLFQ) il nous suffit de 
vérifier l'hypothèse ($\mathcal H2$).
\hfill $\Box$

\nocite{charya}
\nocite{bros}
\nocite{faouzid}
\nocite{faouzic} 
\nocite{f&f&a}
\nocite{duflo}
\nocite{jacod&shirayev}
\nocite{lepingle}
\nocite{lifshits}
\nocite{faiza}
\nocite{revuz}
\nocite{stout}
\nocite{touati1}
\nocite{touati2}
\nocite{wang}
\nocite{damien&gilles}
\bibliographystyle{ims}
\bibliography{biblio}
\end{document}